%% file: multigeodesic.tex
\magnification=\magstep1
\input amstex
\documentstyle{amsppt}

\def\PSL{\operatorname{PSL}}
\def\PSU{\operatorname{PSU}}
\def\SL{\operatorname{SL}}
\def\SU{\operatorname{SU}}

\def\Aut{\operatorname{Aut}}
\font\rsfs=rsfs10

\def\rs#1{\text{\rsfs #1}}

\def\D{\Bbb D}
\def\R{\Bbb R}
\def\Z{\Bbb Z}

\def\ach{\operatorname{arccosh}}

\def\trace{\operatorname{trace}}

\def\ach{\operatorname{arccosh}}
\def\acos{\operatorname{arccos}}
\def\asin{\operatorname{arcsin}}
\def\T{{T\!w}}
\input psfig.sty
\input labelfig.tex

\topmatter
\title
Multi-geodesic tessellations, fractional Dehn twists  and uniformization of 
algebraic curves
\endtitle
\rightheadtext{Multi-geodesic tessellations and uniformization}
\author Samuel Leli\`evre, Robert Silhol
\endauthor
\affil Universit\'e Paris-Sud Orsay, Universit\'e Montpellier II
samuel.lelievre\@math.u-psud.fr, rs\@math.univ-montp2.fr
\endaffil
\address
Universit\'e Montpellier II, D\'epartement de Math\'ematiques,
UMR CNRS 5149,
Place E. Bataillon 34095 Montpellier Cedex 5, France
\endaddress

\endtopmatter
\document
\bigskip
\bigskip

\subhead 0. Introduction \endsubhead 
\medskip

To illustrate what this paper is about we first consider a 
classical example of a translation surface. Consider two copies of a 
regular Euclidean pentagon as on the left of figure 1. Identifying opposite 
parallel sides by translations one obtains a Riemann surface with a natural 
locally flat metric with one cone type point of total angle $6\pi$. Since this 
is a compact Riemann surface it is also an algebraic curve and in this case it 
is well known that this is the curve defined by $y^2=x^5-1$. 
Since the surface is of genus 2 it is also a hyperbolic surface. 
The hyperbolic metric can easily be recovered 
by taking two copies of a regular hyperbolic pentagon (with interior 
angle $\pi/5$) in the unit disk, attached as on the right of figure 1
and with opposite sides identified by appropriate hyperbolic
transformations (this construction is adapted from the one given in 
[Ku-N\"a]).     
\bigskip

\centerline{
 \SetLabels
\endSetLabels
 \AffixLabels{\psfig{figure=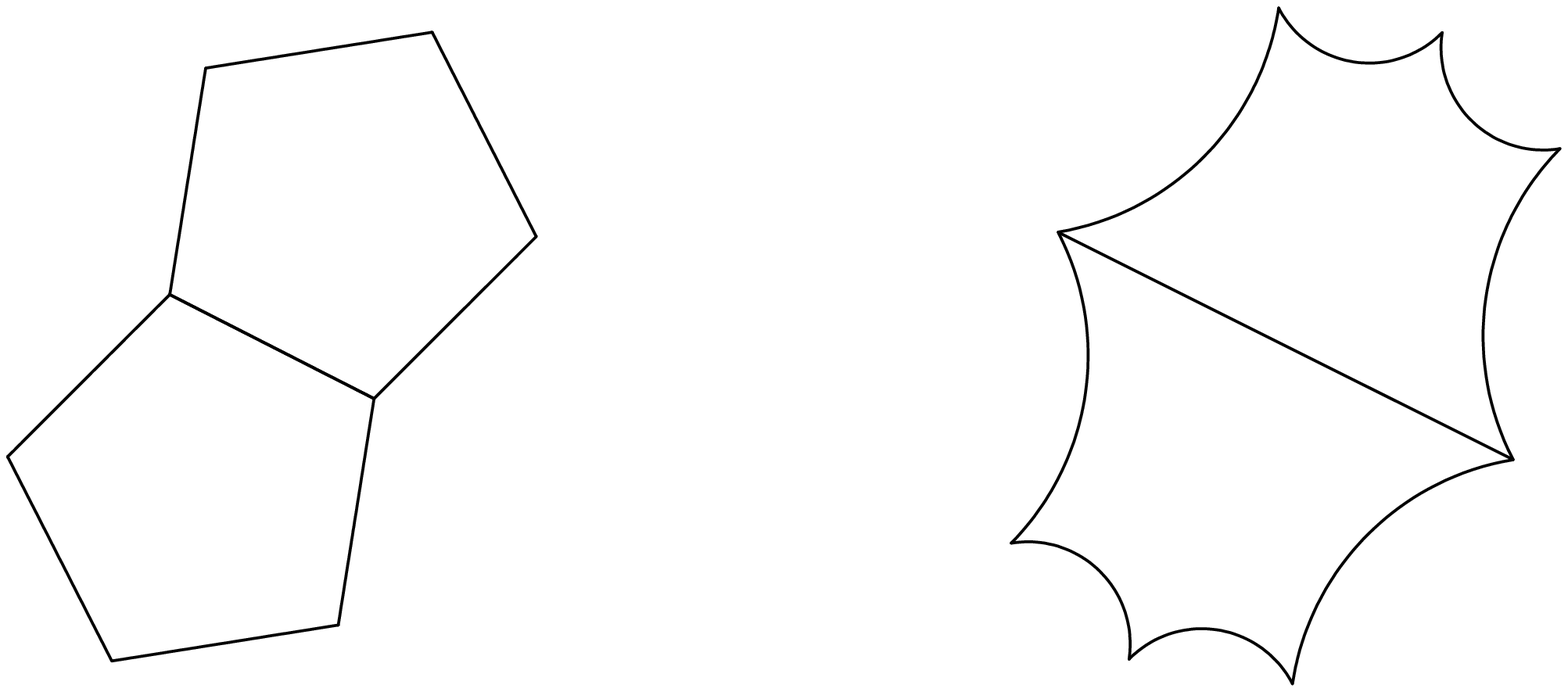,height=3cm}}}
\bigskip
\centerline{\bf Fig. 1}
\bigskip

The important point here is that in fact these two 
decompositions into two regular pentagons coincide exactly or
otherwise said the pentagons are geodesic for both the locally flat 
metric defined above and the natural hyperbolic metric. In this case 
this follows easily from the fact that the edges of the pentagons are
fixed points of anti-conformal reflections. 

This situation is unfortunately far from being the generic one and in 
general two different metrics on a surface have no geodesic arcs in 
common. There are however infinitely many families for which we have a 
decomposition into polygons, geodesic for two or more metrics. This is 
in particular the case for surfaces obtained by taking  finitely many 
copies of an Euclidean rectangle and identifying edges of the 
rectangles pairwise, using translations or rotations of angle $\pi$ 
(half-turns). These are the surfaces we will explore in this paper.
\medskip

The existence of a multi-geodesic tessellation on such surfaces has 
interesting consequences. One of these is that it provides a mechanical 
way to reconstruct a Fuchsian group for the surface (see Propositions 
{\bf 3.2} and {\bf 3.3}).

Another consequence is that it allows for a description in terms of 
Fenchel-Nielsen coordinates of the Teichm\"uller disk generated by a
surface tiled by squares (see Proposition {\bf 3.6}). Moreover such 
a description allows for an interpretation in terms of fractional Dehn 
twists of the natural $\PSL_2(\Bbb Z)$ action on the $\PSL_2(\Bbb R)$-orbit 
of such surfaces, Corollary {\bf 3.7} and section {\bf 5}, where there
are also examples of fractional Dehn twists that connect
surfaces in different $\PSL_2(\Bbb R)$-orbits, or even in different
strata.
These follow from the fact that actually not only is the tessellation 
multi-geodesic but the medians of the rectangles (vertical or horizontal) 
extend to simple closed curves that are geodesic for the different metrics.
\smallskip

Finally we note that 
in many cases one can use the tiling by rectangles to recover an equation 
for the corresponding algebraic curve (see sections {\bf 4} and {\bf 5}). 
Hence the tiling by rectangles provides a bridge between the algebraic 
equation and the  hyperbolic structure deduced 
from the multi-geodesic tessellation. In other words this solves the 
uniformization problem for such curves. In fact it also gives a scheme to do 
uniformization for infinitely many families of curves (see section {\bf 4} 
for some examples). Finally in section {\bf 6} and the Appendix we discuss 
some number theoretic aspects.
\medskip

The authors would like to thank Hugo Akrout and Peter Buser for many useful 
discussions. The authors would also like to thank Le centre  Bernoulli in 
Lausanne for its hospitality while working on this paper.
\bigskip

\subhead 1. Multi-geodesics in genus 1 and genus 0\endsubhead 
\medskip

Let $G_1$ be discrete subgroup of $\Aut(\Bbb D)\cong\PSU(1,1)$ of
genus 1 and generated by two hyperbolic elements $A$ and $B$. The
commutator $[A,B]$ of $A$ and $B$ may be either elliptic of order 
$n$ or parabolic and hence the signature of the group is $(1;n)$ or 
$(1;\infty)$.   

The quotient $S_1=\Bbb D/G_1$ is a hyperbolic genus 1 surface with a cone
point of total angle $2\pi/n$. The surface has also a distinguished
homology basis $(\alpha,\beta)$ given by the images of the (oriented)
axes of the transformations $A$ and $B$. Replacing $B$ by $B^{-1}$ if 
necessary we may assume that this is a canonical basis. In this
situation  there is a unique $\tau$ such that we have a conformal 
equivalence from $S_1$ to $\Bbb C/\Lambda$, where $\Lambda$
is the lattice generated by 2 and $2\tau$, and the image of $[-1,1]$ (resp.\
$[-\tau,\tau]$) under the canonical projection 
$\pi_2:\Bbb C\to\Bbb C/\Lambda=S_1$ is in the same homology class as
$\alpha$ (resp.\ $\beta$). The conditions define the equivalence up 
to a translation in $\Bbb C/\Lambda$. We make it unique by requiring that the 
intersection of the axes maps to $\pi_2(0)$ and write $S_1=\Bbb C/\Lambda$.

Call $\pi_1$ the covering map $\Bbb D\to S_1$. Since $\pi_2:\Bbb C\to S_1$ is 
the universal cover of $S_1$ and $\Bbb D$ is simply connected, $\pi_1$ 
lifts to a map $\varphi:\Bbb D\to\Bbb C$, and this lifting is unique if 
we impose $\varphi(0)=0$.

By construction the surface $S_1$ comes equipped with two metrics, the
natural flat metric and the metric induced by the Poincar\'e metric on
$\Bbb D$. In the general situation the geodesic arcs for these two
metrics bear no relations. A typical situation is represented in
figure 2 where on the left is a fundamental domain in the unit disk
and on the right are the images under $\varphi$ of the hyperbolic
geodesics shown on the left. Although the difference is slight none of
the arcs shown on the right are straight line segments. 
\bigskip

\centerline{
 \SetLabels
\L(.2*.52) $0$\\
\L(.445*.81) $q_1$\\
\L(.03*1.02) $q_2$\\
\L(-.04*.18) $q_3$\\
\L(.37*-.1) $q_4$\\
\L(.39*.5) $p_1$\\
\L(.27*.82) $p_2$\\
\L(.02*.5) $p_3$\\
\L(.15*.16) $p_4$\\
\L(.76*.53) $0$\\
\L(.97*.47) $1$\\
\L(.83*.83) $\tau$\\
\L(1.01*.82) $1+\tau$\\
\endSetLabels
 \AffixLabels{\psfig{figure=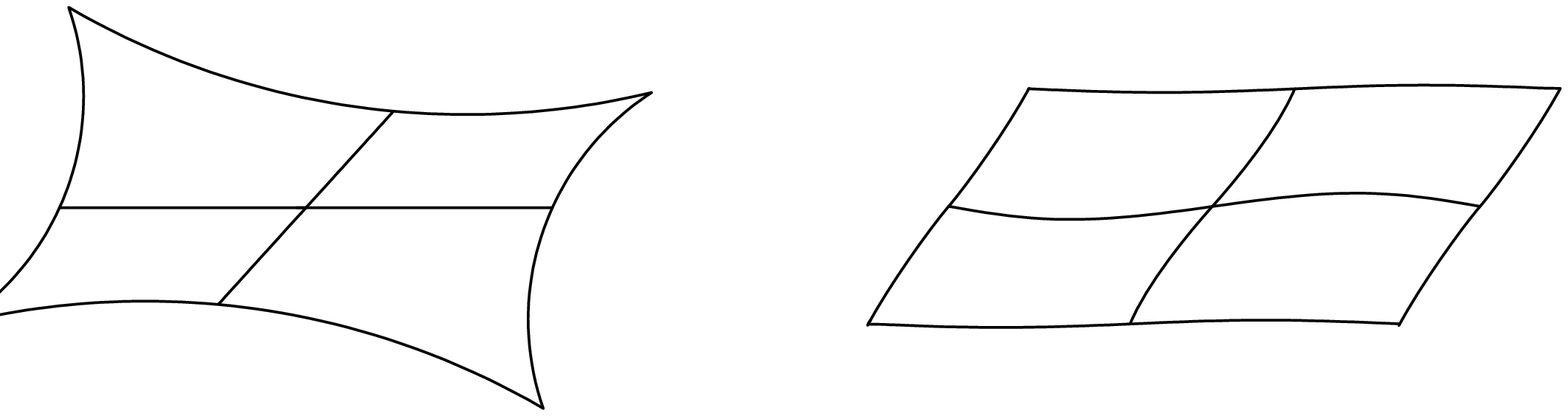,height=2.5cm}}}
\bigskip

\centerline{\bf Fig. 2}
\bigskip

We have nevertheless the following elementary statement 

\proclaim{1.1 Lemma} Let $G_1$, $S_1$, $\tau$ and 
$\varphi:\Bbb D\to\Bbb C$ be as above. Let $q_1=\varphi^{-1}(1+\tau)$,
\dots, $q_4=\varphi^{-1}(1-\tau)$ be the pull-back of the vertices of
the fundamental parallelogram for $S_1$. Then $p_1=\varphi^{-1}(1)$, 
$p_2=\varphi^{-1}(\tau)$,  $p_3=\varphi^{-1}(-1)$ and 
$p_4=\varphi^{-1}(-\tau)$ are the hyperbolic midpoints of $(q_4q_1)$, 
$(q_1q_2)$, $(q_2q_3)$ and $(q_3q_4)$ respectively. Moreover 
$\varphi^{-1}(0)$ is the hyperbolic midpoint of $(p_1p_3)$ and 
$(p_2p_4)$.
\endproclaim

{\smc Proof.} By construction the pull-back via $\varphi$ of $z\mapsto -z$
is the order two elliptic transformation centered at $\varphi^{-1}(0)$ 
and the $q_i$ and $p_i$ are preimages of the fixed points of the induced 
transformation in $S_1$.
\bigskip

If the axes of $A$ and $B$ are orthogonal we have a much stronger statement.  
In this case we may always assume that up to 
conjugation the axis of $A$ is the real axis and the axis of $B$ is the 
pure imaginary axis. The reflections with respect to the 
real axis and the pure imaginary axis induce anti-holomorphic involutions 
on $S_1$ or in other words define real structures on $S_1$. See left of 
figure 3 where a fundamental domain for $G_1$ is represented. In order to 
avoid the denomination ``hyperbolic rectangle'' we will call such a 
domain an {\it equiquadrangle}.
\bigskip
\medskip

\centerline{
 \SetLabels
\L(.12*.4) $0$\\
\L(.27*.49) $p_1$\\
\L(.31*1.02) $q_1$\\
\L(.14*.94) $p_2$\\
\L(.85*.4) $0$\\
\L(1.01*.48) $1$\\
\L(1.01*.86) $1+\tau$\\
\L(.875*.9) $\tau$ \\
\endSetLabels
 \AffixLabels{\psfig{figure=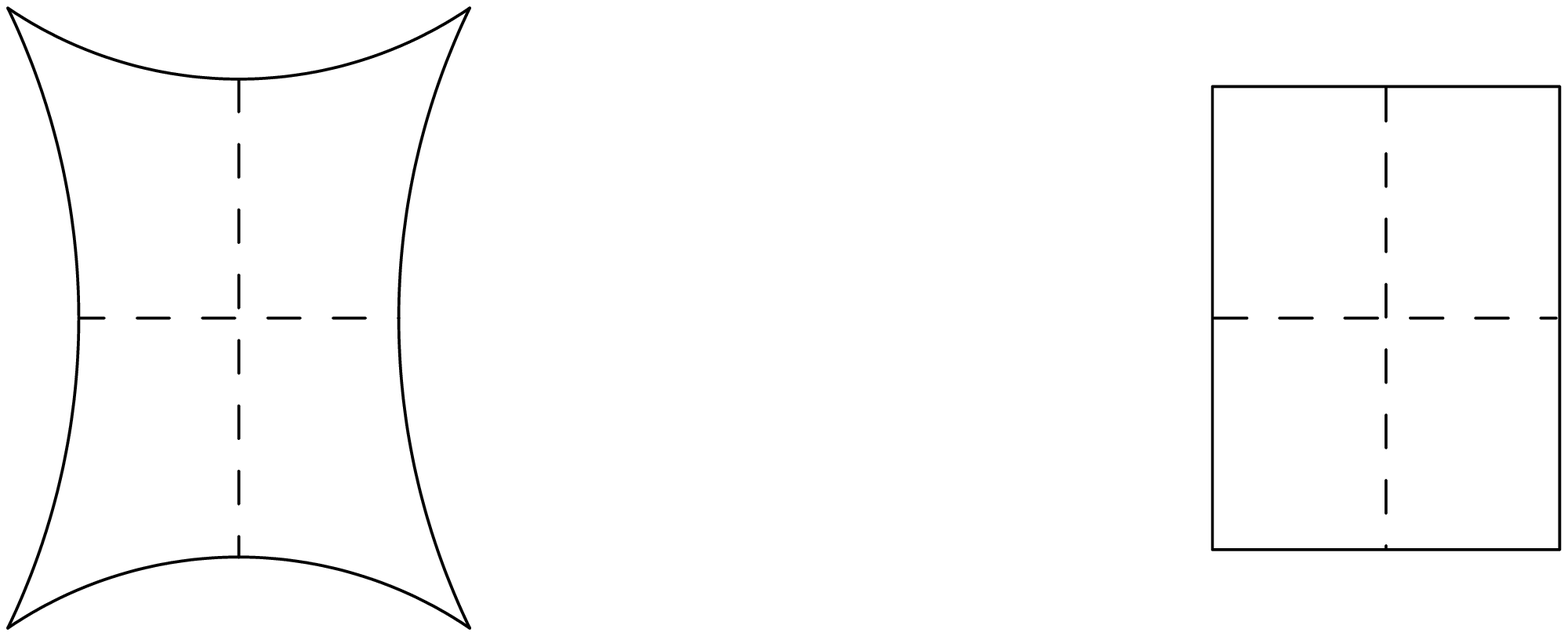,height=3cm}}}
\bigskip
\centerline{\bf Fig. 3}
\bigskip

These real structures obviously have two real components, one is the image 
of the real (resp.\ pure imaginary) axis, the other is the image of two 
identified opposite sides of the fundamental domain. Looking at the
lattice $\Lambda=\langle 2,2\tau\rangle$ introduced above, this means that
$\tau$ is pure imaginary, or in other words that the natural fundamental 
domain for $\Lambda$ is a rectangle. By the uniqueness of the map
$\varphi$, this map commutes  with complex conjugation on both sides, 
i.e.~ $\overline{\varphi(\overline{z})}=\varphi(z)$ and hence maps 
$\Bbb D\cap\Bbb R$ to $\Bbb R$ and $\Bbb D\cap i\,\Bbb R$ to $i\,\Bbb R$. 

\proclaim{1.2 Lemma}  Let $G_1$ be a discrete subgroup of $Aut(\Bbb D)$ 
of signature $(1;n)$ or $(1;\infty)$ and generated by two hyperbolic
elements $A$ and $B$ with respective axes the real axis and the pure
imaginary axis. Let $A_1$ (resp.\ $B_1$) be the unique hyperbolic
element such that ${A_1}^2=A$ (resp.\ ${B_1}^2=B$). 

Let $c_n$ be the images of the pure imaginary axis under ${A_1}^n$ 
and $d_n$ the image of the real axis under ${B_1}^n$. 
Let $S_1=\Bbb D/G_1$ and $\pi_1$ the natural projection $\Bbb D\to S_1$. 
Then for all $n$, $\pi_1(c_n)$ and $\pi_1(d_n)$ are geodesic
arcs in $S_1$ for both the hyperbolic metric induced by that of $\Bbb D$ 
and for the natural conformal flat metric on the torus $S_1$.
\endproclaim  
  
\noindent{\smc Proof:} The real structure on $S_1$ induced by complex
conjugation in $\Bbb D$ has two real connected components 
$\pi_1(\Bbb R\cap\Bbb D)$ and $\pi_1(d_1)$. Since these are the fixed 
points of an anti-conformal reflection they are geodesic for any
metric compatible with the conformal structure and for which the
reflection is anti-conformal. The same is true for 
$\pi_1(i\Bbb R\cap\Bbb D)$ and $\pi_1(c_1)$. Since the $c_n$ and $d_n$ 
are just the images of $c_0$, $c_1$, $d_0$ and $d_1$ under $G_1$ we are 
done.
\bigskip

\proclaim{1.3 Corollary} Let $S_1$ be as in {\bf 1.2} and let 
$\varphi:\Bbb D\to\Bbb C$ be the lifting of the projection map 
$\pi_1:\Bbb D\to S_1$ normalized as above. Let $G$ be the group
generated by $A_1$ and $B_1$ --- see  {\bf 1.2}. Then the images 
of $\Bbb D\cap\Bbb R$ and $\Bbb D\cap i\,\Bbb R$ under $G$ are mapped
by $\varphi$ onto vertical lines through the integers and horizontal 
lines through $n\,\tau$, $n\in\Bbb Z$.
\endproclaim

Another case where we have multi-geodesic arcs is the case of genus 1 
surfaces with a half-Dehn twist (see figure 4).

\proclaim{1.4 Lemma} Let $A$, $B$, $A_1$, $B_1$ and $G_1$ be as in Lemma 
{\bf 1.2}. 
Let $G'_1$ be the group generated by $A$ and $B'={A_1}\,B$. Then 
$G'_1$ has same signature as $G_1$ (i.e.~ $(1;n)$ or $(1;\infty)$). Let 
$S'_1=\Bbb D/G'_1$ and let $\pi'_1:\Bbb D\to S'_1$ be the natural projection.

Then $\pi'_1(i\Bbb R\cap\Bbb D)$, $\pi'_1({A_1}(i\Bbb R\cap\Bbb D))$ and 
 $\pi({B_1}^{-1}(\Bbb R\cap\Bbb D))$ are geodesic arcs in $S'_1$ for both the
hyperbolic metric induced by $\Bbb D$ and the natural flat metric on the
torus.
\endproclaim

\noindent{\smc Proof.} Since $A$ and $A_1$ commute, the commutator $[A,B']$
is conjugate to the commutator $[A,B]$. This proves the first assertion.

For the rest we note that $S'_1$ is obtained from $S_1$ by applying a
half-Dehn twist along the image of the axis of $A$. Since this is a real 
component of $S_1$, $S'_1$ also has real structures, but with only one 
real component this time (see for example [Bu-Se]). The real structures 
compatible with the hyperbolic metric must keep fixed the elliptic point 
(or the cusp, depending on the signature). To describe these, let $\sigma$ 
be complex conjugation and let $\sigma_1=B_1^{-1}\cdot\sigma\cdot B_1$ and 
$\sigma_2=-\sigma$. We have $B_1^{-1}\cdot\sigma=\sigma\cdot B_1$, hence we 
have $\sigma_1=B^{-1}\cdot\sigma=\sigma\cdot B$. We also have 
$\sigma\cdot A=A\cdot\sigma$ and similarly $\sigma\cdot A_1=A_1\cdot\sigma$.
With these relations it is easy to prove that
\roster
\item"(i)" $\sigma_1\cdot A\cdot\sigma_1={B'}^{-1}\cdot A\cdot B'$;
\item"(ii)" $\sigma_1\cdot B'\cdot\sigma_1={B'}^{-1}\cdot A$.
\endroster
In exactly the same way we can also prove
\roster
\item"(iii)" $\sigma_2\cdot A\cdot\sigma_2=A^{-1}$;
\item"(iv)" $\sigma_2\cdot B'\cdot\sigma_2={A}^{-1}\cdot B'$.
\endroster
 
To end the proof we only need to note that the fixed part of $\sigma_1$ 
is ${B_1}^{-1}(\Bbb R\cap\Bbb D)$ and the fixed part of $\sigma_2$ is 
$i\Bbb R\cap\Bbb D$ which has same image as ${A_1}(i\Bbb R\cap\Bbb D)$.
\bigskip
\medskip

\centerline{
 \SetLabels
\endSetLabels
 \AffixLabels{\psfig{figure=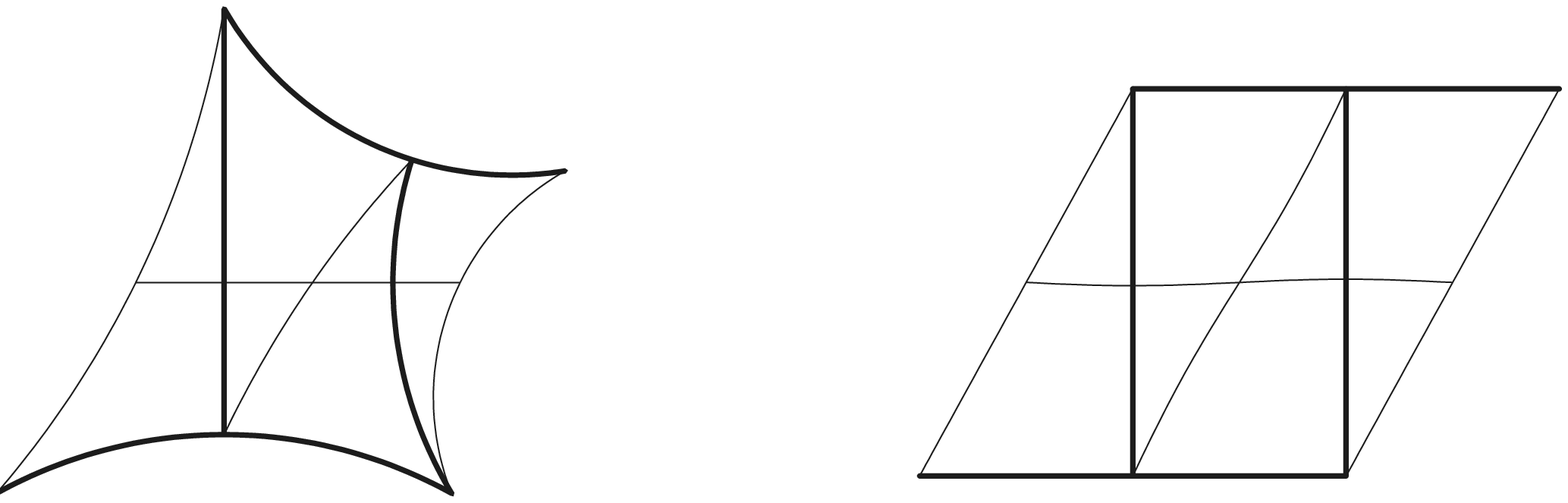,height=3cm}}}
\bigskip
\centerline{\bf Fig. 4}
\bigskip
\bigskip
 
In general no other geodesic gets mapped in $\Bbb C$ onto a straight line. 
There are however some noteworthy exceptions. Of particular interest are 
the quadrangles corresponding to $\tau=i$, $\tau=\frac12+\frac{i}2$ and 
$\tau=(1+i\,\sqrt{3})/2$. In these cases the existence of additional
 real structures implies that the arcs marked  in figure 5 are also 
geodesic for the hyperbolic metric.
\bigskip

\centerline{{\psfig{figure=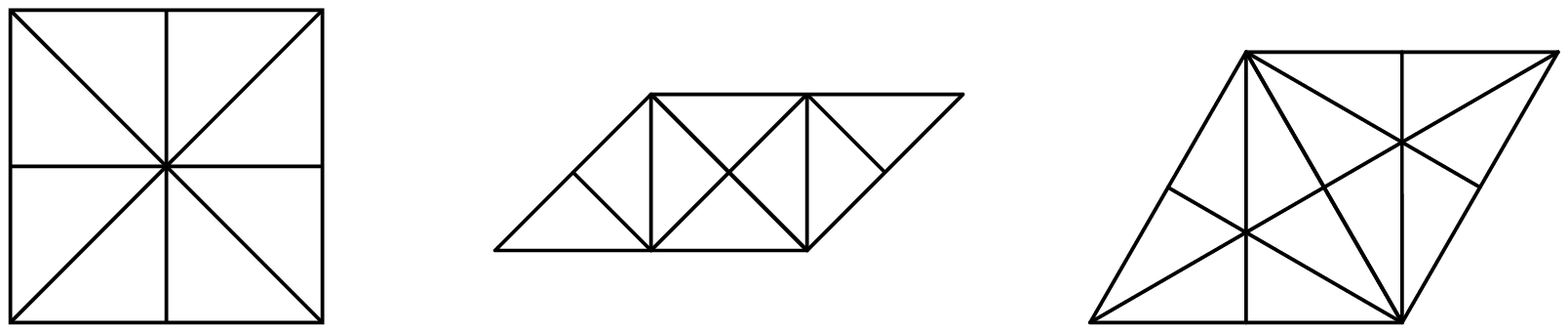,height=2cm}}}
\bigskip

\centerline{\bf Fig. 5}
\bigskip

In the sequel we will need to use a genus 0 variant of {\bf 1.2}.
Let $A$ be as above a hyperbolic element with axis the real line. 
Let $e_1$ be the elliptic element of order 2 with center 0 and let 
$e_2$ be a second elliptic element of order 2 with center in a point
$\tau$ on the pure imaginary axis. This data being subject to the condition
that the group $G_0=\langle A,e_1,e_2\rangle$ is a discrete subgroup 
of $\Aut(\Bbb D)$ with signature of the form $(0;2,2,2,n)$ or 
$(0;2,2,2,\infty)$. By construction the quotient $\Bbb D/G_0$ is of 
genus 0 and is naturally equipped with a singular hyperbolic metric with 
4 cone points, 3 with total angle $\pi$ and one with angle $2\pi/n$ or 
3 cone points of angle $\pi$ and a cusp.
\bigskip

\centerline{
 \SetLabels
\L(.16*.225) $0$\\
\L(.26*0.1) $p_1$\\
\L(.31*.86) $q_1$\\
\L(.14*.76) $p_2$\\
\L(.9*.22) $0$\\
\L(1.01*0.15) $1$\\
\L(1.01*.82) $1+\tau$\\
\L(.88*.85) $\tau$\\
\endSetLabels
 \AffixLabels{\psfig{figure=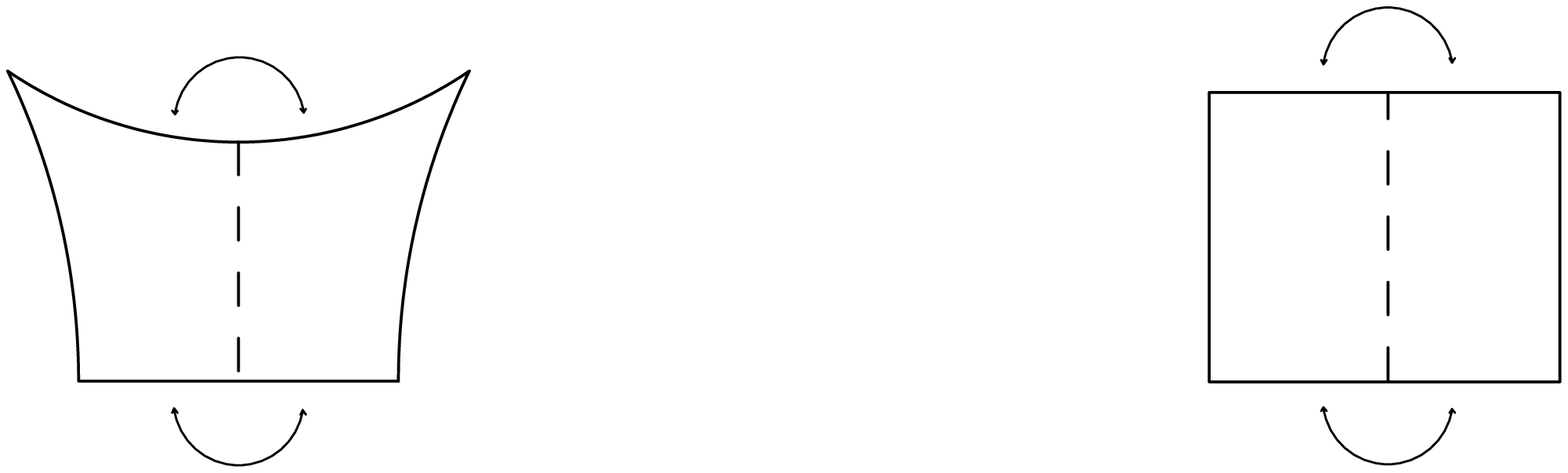,height=2.5cm}}}
\bigskip

\centerline{\bf Fig. 6}
\bigskip

The fundamental domain described on the left of figure 6 is conformally
equivalent to a rectangle with vertices $\pm 1$ and $\pm 1+\tau$. 
Using this rectangle we can also equip $S_0$ 
with a singular flat metric (see right of figure 6). For the same reasons 
as those indicated for $S_1$ the sides of the fundamental domain coincide
with the sides of the rectangle, and in particular are geodesic for both 
metrics. For further use we write this formally.

\proclaim{1.5 Lemma} Let $G_0$ be a discrete subgroup of $\Aut(\D)$ of 
signature $(0;2,2,2,n)$ or $(0;2,2,2,\infty)$ and generated by an 
elliptic element $e_1$ of order $2$ centered at 0 and two hyperbolic
elements $A$ and $B$ with respective axes the real axis and the pure
imaginary axis.

Let $S_0=\D/G_0$. Then $S_0$ decomposes into 
two copies of a geodesic hyperbolic trirectangular quadrangle. 
The sides of this quadrangle are also geodesic for a natural singular flat 
metric.
\endproclaim

We have just taken $B=e_1e_2$.
\bigskip 

\noindent{\bf 1.6 Remark.} It will be useful in the sequel to reformulate
the conditions on $G_1$ and $G_0$ in terms of fundamental domains. 
For $G_1$ the conditions are that a fundamental domain in the disk 
be of the form given on the left of figure 3, with interior angles 
all equal either to $\pi/(2n)$ or to 0. For $G_0$ the condition is 
that the fundamental domain be of the form given on the left of 
figure 6, with two interior right angles and two angles equal either
to $\pi/n$ or to 0.
\bigskip 

\subhead 2. Fuchsian groups and equations in genus 1\endsubhead
\medskip

For the applications we will need to have a more precise description of 
the hyperbolic transformations $A$ and $B$ introduced in Lemma {\bf 1.2}
and their relations with the rectangles.
\smallskip

\proclaim{2.1 Lemma} Let $\rs R$ be an equiquadrangle with interior
angle $\pi/n$, $n$ even, or zero angle.
Let $\ell$ be the hyperbolic length of the horizontal median of $\rs R$ 
and let $L=\cosh(\ell/2)$. Let 
$$L'=\sqrt{\frac{cos(\pi/n)^2+L^2-1}{L^2-1}}\ \text{ or }\ 
L'=\frac{L}{\sqrt{L^2-1}}\ \text{ if the angle is zero.}$$

Then the transformations $A$, $B$ of Lemma {\bf 1.2} generating $G_1$
are represented in $\SU(1,1)$ by
$$A=\pmatrix L & \sqrt{L^2-1}\\ 
\sqrt{L^2-1} & L\endpmatrix,\qquad 
B=\pmatrix L' & i\,\sqrt{{L'}^2-1}\\ -i\,\sqrt{{L'}^2-1} & L'\endpmatrix
\ .\tag 2.1.1$$

In the more general situation of genus 1 with a twist parameter $t$ along 
the axis of $A$ we can consider the group generated by $A$ and 
$B_1=T\cdot B$, where
$$T=\pmatrix Tw & \sqrt{Tw^2-1} \\ \sqrt{Tw^2-1} & Tw \endpmatrix\tag 2.1.2$$
with $Tw=\cosh(t\,\ach(L))$. 
\endproclaim

The proof is a simple exercise in hyperbolic trigonometry (use for example 
[Bu], p.454). 
\smallskip

For the group $G_0$ of Lemma {\bf 1.5} we have a very similar description.

\proclaim{2.2 Lemma} Let $G_0=\langle A,B,e_1\rangle$ be as in Lemma 
{\bf 1.5}. Then $A$ and $B$ will have the same expression as that given in 
$(2.1.1)$, but with $n$ even or odd this time. 

In the presence of a twist $t$ along the axis of $A$, $A$ and $B_1$ have the 
same expression as in Lemma {\bf 2.1} and $e_1$ is just the conjugate of 
$z\mapsto -z$ by a matrix of the same form as $T$ but with $Tw$ replaced by 
$\sqrt{(Tw+1)/2}$.
\endproclaim 

For the special cases of quadrangles associated to the parallelograms of 
periods $i$, $\frac12+\frac{i}2$ or $(1+i\,\sqrt{3})/2$, we can be even 
more precise. 
\smallskip

\noindent{\bf (2.3)} For the ``square'' equiquadrangle we only need to take
$$L=\sqrt{\cos(\pi/n)+1}\ \text{ or }\ L=\sqrt{2}\ \text{  
if the angle is 0.}\tag 2.3.1$$

\noindent{\bf (2.4)} If $\tau=(1+i\,\sqrt{3})/2$ the corresponding hyperbolic
   quadrangle is obtained from two copies of a hyperbolic equilateral
   triangle with angle $2\pi/(3n)$. Using the formulae in [Bu], p.454, 
   one can compute the generators $A$ and $B$ for the group of the
   genus 1 surface. This yields:
   \roster
   \item"$\bullet$" $A$ as in (2.1.1) with 
   $L=\frac12+\cos\left(\frac{2\pi}{3n}\right)$;
   \item"$\bullet$" $B=T\,A\,T^{-1}$,
   \endroster
   where 
   $$T=\pmatrix \exp\left(i\,\frac{\theta}2\right) & 0\\ 0 & 
   \exp\left(-i\,\frac{\theta}2\right)\endpmatrix
   \quad \text{ with }\ 
   \theta=\acos\left(\frac{2\cos(2\pi/(3n))+1}{2\cos(2\pi/(3n))+3}\right)\ .$$

   If the angle is the zero angle these matrices become
   $$A=\pmatrix \frac32 & \frac{\sqrt{5}}2\\ \frac{\sqrt{5}}2 &\frac32
   \endpmatrix\ \text{ and }\ B=\pmatrix \frac32 & 
   \frac{(3+4\,i)\sqrt{5}}{10}\\ \frac{(3-4\,i)\sqrt{5}}{10} &\frac32
   \endpmatrix\ .$$

\noindent{\bf (2.5)} If $\tau=\frac12+\frac{i}2$, then one has, with $L$ as in 
(2.3.1),
$$A=\pmatrix L^2& \sqrt{L^4-1}\\ \sqrt{L^4-1} & L^2\endpmatrix,\quad  
B=\pmatrix L &e^{i\beta}\,\sqrt{L^2-1}\\ e^{-i\beta}\,\sqrt{L^2-1}& L
\endpmatrix$$
where $\beta=\asin(1/\sqrt{L^2+1})$.
\bigskip

We will also need to associate an explicit equation for 
the elliptic curve defined by a rectangle and conversely.

We do this as 
follows, given $\tau$ we use a variant of a classical Jacobi function
$$J_{\tau}(z)=-\frac{w}{K}\prod_{k=0}^{\infty}
\frac{(w-\zeta^{2k})^2(1-\zeta^{2k+2}w)^2}
{(w-\zeta^{2k+1})^2(1-\zeta^{2k+1}w)^2}\ ,\tag{2.6}$$
where $\zeta=\exp(\pi i \tau)$, $w=\exp(\pi i z)$ and 
$K=4\prod\limits_{k=1}^{\infty}
\left(\frac{1+\zeta^{2k}}{1+\zeta^{2k-1}}\right)^4$
(for more details see for example Nehari [Ne], Chap. VI, Sec.3. See also 
[Bu-Si2]). 
We have $J_{\tau}(0)=0$,  $J_{\tau}(1)=1$ and  $J_{\tau}(\tau)=\infty$. 
Letting $\mu=J_{\tau}(1+\tau)$, an equation for the elliptic curve is
$$y^2=P(x)=x\,(x-1)\,(x-\mu)\ .\tag 2.7$$

For the converse choose a determination of $\sqrt{P(x)}$ for $x$ in the 
upper half plane. Extend this determination to the real line. Then
$$\tau=I_2/I_1\ \text{ where }\ I_1=\int_0^1\frac{dx}{\sqrt{P(x)}}\ 
\text{ and }\ I_2=\int_1^{\mu}\frac{dx}{\sqrt{P(x)}}\ .$$
\medskip

We have used here the convention that the 
vertices of the rectangle get mapped to the Weierstrass point with 
$x$-coordinate 0, the midpoint of the vertical edges to the one with 
$x$-coordinate 1 and the midpoint of the horizontal edges to the point 
at infinity. We extend this convention to the case of parallelograms 
namely the vertices will be mapped to 0, the midpoint of the horizontal 
edges to $\infty$ and the midpoints of the other edges to 1.

With this convention we note that
\proclaim{2.8 Lemma} If $\mu>1$ corresponds to $\tau$, then 
$\dfrac{\mu}{\mu-1}$ corresponds to $-1/\tau$, $1-\mu$ corresponds to 
$\tau-1$ and $1/\mu$ to $\dfrac{\tau}{1+\tau}$.
\endproclaim

\noindent{\smc Proof.} For the first assertion we note that the integrals 
along $[-\infty,0]$ (resp. $[\mu,\infty]$) are just minus the integrals 
along $[1,\mu]$ (resp. $[0,1]$). Applying the change of variable 
$x\mapsto x/(x-1)$ proves the first assertion if in addition we note that 
the integral along $[0,1]$ is real while the integral along $[1,\mu]$ is 
pure imaginary. For the other assertions apply the change of 
variables $x\mapsto 1-x$ and $x\mapsto 1/x$.
\bigskip 

\noindent{\bf 2.9 Remarks.} 1) If we are dealing with a parallelogram in 
place of a rectangle, then the curve will also have an equation of the 
form (2.7) with $\mu$ computed using the Jacobi function (2.6). We will 
also of course have an analogue of Lemma {\bf 2.8} but the exact 
correspondence will depend on the choice of the initial parallelogram.

2) Since for a given elliptic curve there 
are only six possibly different values of $\mu$, infinitely many values 
of $\tau$ will correspond to the same $\mu$. This is in particular the
case for $\tau$ and $2+\tau$.
\medskip 

3) For domains as in the second part of Lemma {\bf 2.1} 
we use a similar convention that 0 is mapped to $0$, 1 is mapped 
to 1 and $\tau$ is mapped to $\infty$. We will call $\mu$ the image of 
$1+\tau$.
\bigskip

A fundamental domain for a group generated  $A$ and $B$ (resp. $A$ and 
$T\cdot B$) is illustrated on the left (resp. right) of figure 7.
  \bigskip
  \bigskip

\centerline{
 \SetLabels
\L(-.03*-.07) $q_3$\\
\L(.3*-.07) $q_4$\\
\L(-.03*.98) $q_2$\\
\L(.3*.98) $q_1$\\
\L(.14*0) $p_4$\\
\L(.14*.91) $p_2$\\
\L(-.02*.45) $p_3$ \\
\L(.28*.45) $p_1$\\
\L(.6*-.07) $q_3$\\ 
\L(.76*0) $p_4$\\
\L(.92*-.07) $q_4$\\
\L(.7*1.05) $q'_2$ \\
\L(.86*.89) $p'_2$ \\
\L(1.01*.87) $q'_1$\\
\endSetLabels
 \AffixLabels{\psfig{figure=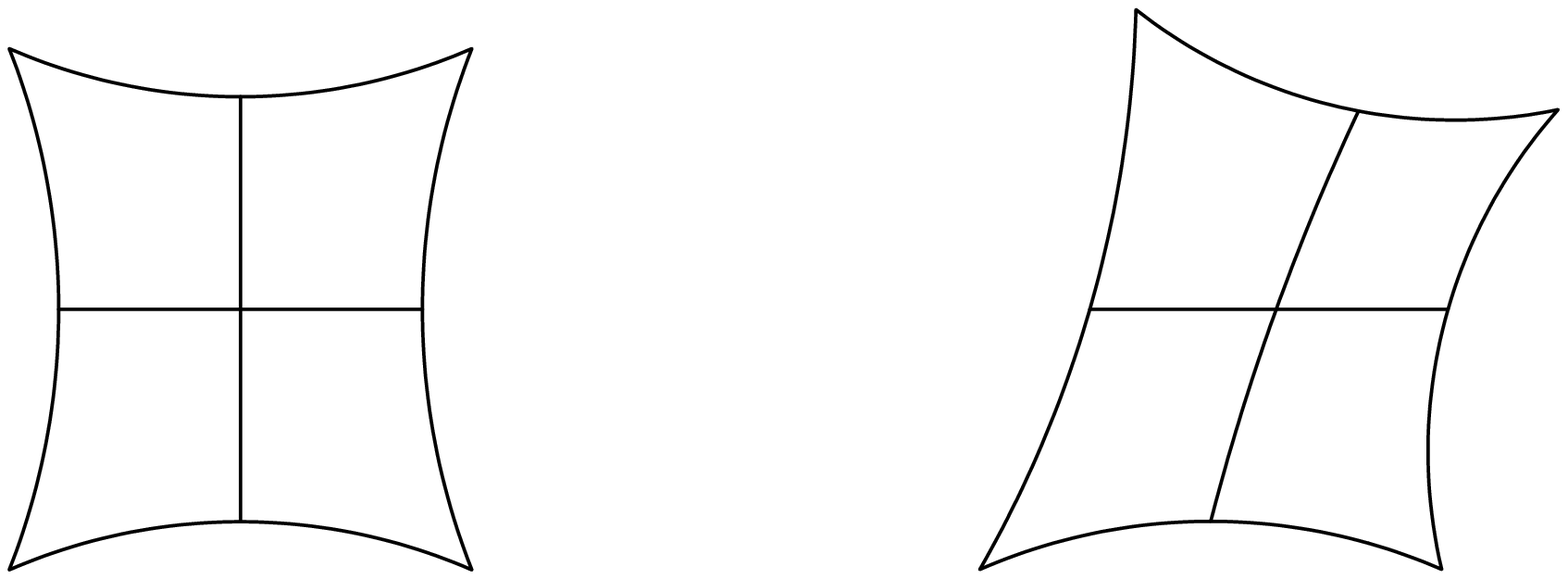,height=3cm}}}
\bigskip
\centerline{\bf Fig. 7}
\bigskip
 
We will need the following straightforward result

\proclaim{2.10 Lemma} Let $L$, $n$, $A$, $B$ and $T$ be as in Lemma 
{\bf 2.1} or Lemma {\bf 2.2}. Then the hyperbolic lengths of the upper 
and lower geodesic arcs of the fundamental domains for 
$\langle A,B\rangle$ or $\langle A,T\cdot B\rangle$ is 
$\ach(L\,\sin(\pi/n))$. Moreover the sum of the interior angles at
$q'_1$ and $q'_2$ is $2\,\pi/n$.
\endproclaim

{\smc Proof.} This follows immediately from the fact that $q'_1$ 
and $q'_2$ are the images of $q_1$ and $q_2$ under $T$.

\bigskip

For the ``square'' equiquadrangle or the quadrangles corresponding to 
the period $\frac12+\frac{i}2$ or $(1+i\,\sqrt{3})/2$ we have the obvious  
values $\mu=2$, $\mu=1/2$ and $\mu=(1+i\,\sqrt{3})/2$ respectively.

In addition to these values there are other cases for which we can 
express $\mu$ in terms of $L$ and $n$. For some examples see the 
Appendix.
\bigskip

\subhead 3. Multi-geodesic tessellation of surfaces obtained from 
rectangles\endsubhead
\medskip

Let $R$ be an Euclidean rectangle in the complex plane. To
simplify we assume the edges of $R$ are horizontal or vertical. 
Assemble $r$ copies $R_1,\dots, R_r$ of $R$ by pasting
along sides of same length so that not only the resulting polygon is
simply connected but it also remains simply connected when one removes
the vertices of the rectangles. We will say that such an arrangement
of rectangles is an {\it admissible arrangement} (the arrangement on 
the left of figure 8 is not admissible in this sense but the one on the 
right is). This restriction on the arrangement is here for purely
technical reasons and could be in fact dispensed with.
\bigskip

\centerline{
 \SetLabels
\endSetLabels
 \AffixLabels{\psfig{figure=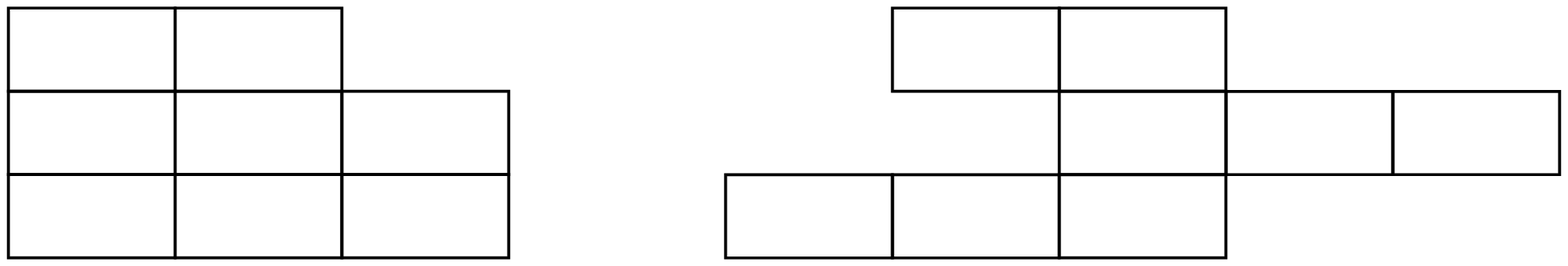,height=1.8cm}}}

\hskip 1cm non admissible\hskip 4cm admissible
\bigskip

\centerline{\bf Fig. 8}
\bigskip
 
Identify the remaining edges by pairs, using translations or rotations 
of angle $\pi$ (half-turns), in a way that is compatible 
with the orientation of the rectangles, i.e.~ the identifications by
translation can be top edge to bottom edge and right edge to left
edge, and the identifications by half-turns can be top
edge to top edge, left edge to left edge, right edge to right edge and
bottom edge to bottom edge.
We will call such surfaces {\it surfaces obtained from rectangles}.  

We will distinguish the case when the identifications are 
only by translations and the case when some identifications are by 
half-turns. Referring to the underlying
invariant differential we will call the first the {\it Abelian case} 
and the second the {\it quadratic case}. This denomination comes from
the fact that in the first case the differential $dz$ being
invariant by translations it induces a holomorphic, or abelian, differential
$\omega$ on the surface $S$. In the second it is $dz^2$ that induces a
quadratic differential $q$. In both cases the zeros of 
the differential are necessarily located at images of vertices of the 
rectangles. 

\proclaim{3.1 Definition} Let $S$ be as above and let $q$ be the
quadratic differential induced by $dz^2$ (in both the Abelian and
quadratic case). 
Let $d_1,\dots,d_k$ be the orders of the zeros of $q$. Let 
$m_i=d_i+2$ and let $m$ be the least common multiple of the $m_i$. 
For a positive integer $n$ we will say that $\pi/n$ (or by abuse 
simply $n$) satisfies the angle condition for the differential if $n$ 
is a multiple of $m$.  We will always consider that the zero angle 
satisfies the angle condition. 
\endproclaim

We consider first the Abelian case.
Let $E$ be the genus 1 surface corresponding to the rectangle $R$ 
that to fix notations we identify with $R_1$ and let $E^*$ be the
surface obtained from $E$ by removing the image of the vertices of the
rectangle.  Let $p$ be the center of $R_1$ and let $h$ and $v$ be 
the elements of $\pi_1(E^*,p)$, 
the fundamental group of $E^*$, corresponding respectively to the
horizontal median and the vertical median of the rectangle oriented in the 
natural way. Since $\pi_1(E^*,p)$ is isomorphic to the free group in
two generators, $h$ and $v$ can be thought of as generating a free
group. 

Let $S$ be as above and let $S^*$ be the surface obtained by removing
the vertices of the rectangles. Since the arrangement of rectangles we
are starting with is admissible and in particular remains simply
connected when the vertices are removed the identifications of edges
define a set of generators of the fundamental group of $S^*$. 
On the other hand $S^*$ is an unramified covering of $E^*$. Hence we
can consider the fundamental group of $S^*$ as a subgroup of
$\pi_1(E^*,p)$. More precisely taking as base point the center of the 
rectangle $R_1$ the fundamental group of $S^*$ is generated by
words $w_1(h,v),\dots w_s(h,v)$ in $h$ and $v$. For example if we 
consider the surface obtained from three rectangles (see upper left of 
figure 9) with the usual identifications (by horizontal and vertical 
translations) then these words are 
$h^2,v^2,hvh^{-1},vhv^{-1}$ (see [Sch] and section {\bf 4} for other 
explicit computations of this type).
\medskip

Let $S$ and $\omega$ be as above and let $n$ be an integer satisfying
the angle condition for $\omega^2$. Then there exist hyperbolic
elements $A$ and $B$ satisfying the conditions of {\bf 1.2} generating
a group $G_1$ of signature $(1;n)$ or $(1;\infty)$ and such that 
$\Bbb D/G_1=E$.
\bigskip

\centerline{
 \SetLabels
\endSetLabels
 \AffixLabels{\psfig{figure=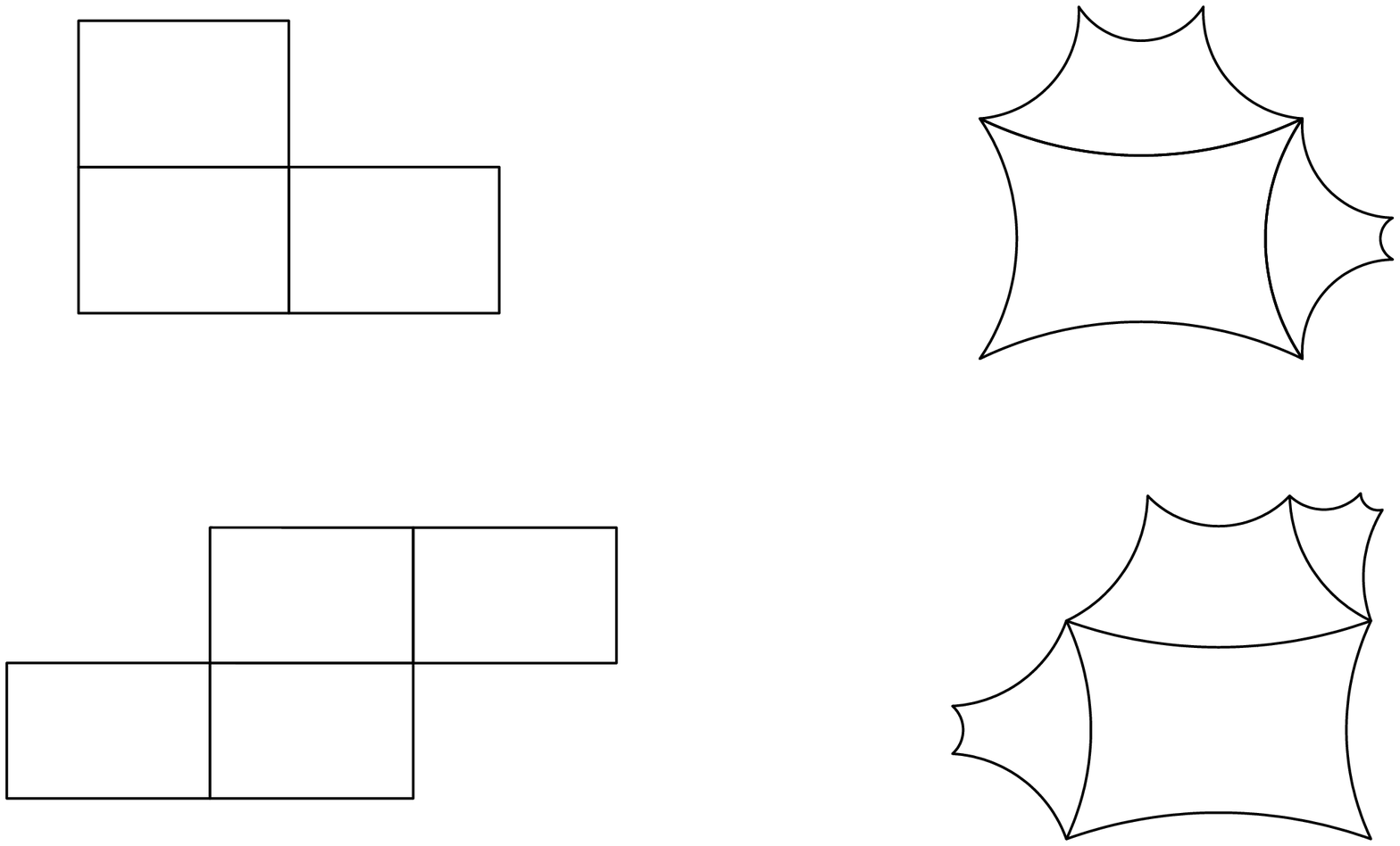,height=4cm}}}
\bigskip

\centerline{\bf Fig. 9}
\bigskip

\proclaim{3.2 Proposition} Let $S$, the Abelian differential $\omega$,
the words $w_1,\dots,w_s$ and $A$ and
$B$ be as above. Let $G$ be the group generated by
$w_1(A,B),\dots,w_s(A,B)$. Then
\roster
\item"(i)" $G$ is a discrete subgroup of $\Aut(\Bbb D)$;
\item"(ii)" $\Bbb D/G\cong S$ as Riemann surfaces;
\item"(iii)" the edges of the rectangles are geodesic arcs for the
hyperbolic metric induced by the Poincar\'e metric on $\Bbb D$;
\item"(iv)" the horizontal and vertical medians of each rectangle are
geodesic for the hyperbolic metric and extend to simple closed
geodesics on the surface.
\endroster
\endproclaim

\noindent{\smc Proof.} Call $\rs R$ the equiquadrangle
defined by $A$ and $B$. We build a domain in $\D$ starting with 
$\rs R$ and using the same combinatorics as the one that defines the
arrangements of rectangles, replacing everywhere the copies of the Euclidean
rectangle $R$ by copies of $\rs R$  (see figure 9 for two examples). By
construction $G$ identifies in pairs the remaining edges of the
equiquadrangles. 
  
Assume that the signature of $G_1=\langle A,B\rangle$ is $(1;n)$.
If $x_0$ is a point of $S$ where $\omega$ has a zero of order $d_i$,  
the total angle at that point will be $2(d_i+1)\pi$ and hence
corresponds to the identification of $4(d_i+1)$ vertices.
The interior angles of the equiquadrangle are $\pi/(2n)$ hence
the sum of the hyperbolic interior angles of the vertices identified
by $G$ with $x_0$ will be $4(d_i+1)\pi/(2n)$. Since by
definition $n$ is a multiple of $d_i+1$ this is of the form $2\pi/p$.
This together with the pairing convention means that
we can apply  Poincar\'e's theorem (see for example
[Bea] section 9.8) and conclude that the group $G$ is discrete and
that the arrangement of equiquadrangles we have constructed is a
fundamental domain for this group. If the signature is $(1;\infty)$ 
then we can directly apply Poincar\'e's theorem (see [Bea], p.251).
This proves (i). One should think of $G$ as an orbifold fundamental
group. 

Next we consider the conformal equivalence $f$ from the interior
of the equiquadrangle to the rectangle, extended by continuity to the 
boundary minus the vertices. Using Schwarz's reflection principle we can
extend this conformal map to the full arrangement of equiquadrangles.
This yields a conformal equivalence from the hyperbolic arrangement to
the Euclidean arrangement. By construction of $G$ this conformal
equivalence induces a conformal equivalence from $\D/G-\{$images of
vertices$\}$ to $S^*$ which extends naturally to a biholomorphic map 
from $\D/G$ to $S$. This proves (ii) and (iii).  

The median (horizontal or vertical) of a rectangle in $S$ extends to a
simple closed curve which corresponds to the decomposition of $S$ into
cylinders (horizontal or vertical). By {\bf 1.2} these are locally
geodesic and since they intersect orthogonally at midpoints the sides of the
equiquadrangles (see {\bf 1.1}) we obtain (iv).     
\medskip

In the quadratic case we have a very similar statement. Essentially
only the initial setup is different.

By changing the arrangement if necessary we may always assume that the
identifications of the form $z\mapsto -z+c$ are between horizontal
sides. 

Let $h$ be the translation that maps the left side of the rectangle
$R$  onto the right side. Let $r_1$ be the rotation of angle $\pi$ 
centered at the center of $R$ and let $r_2$ be the rotation of 
angle $\pi$ centered at the middle of the upper side of $R$. Let 
$S_0^*$ be the quotient of $R-\{$vertices$\}$ under the 
identification induced by $h$, $r_1$ and $r_2$. We consider $S_0^*$ as 
an orbifold with three cone points of order 2 and a cusp. The orbifold 
fundamental group of $S_0^*$ is generated by three elements, two of order 
2 and one of infinite order. We may consider these generators as being  
$r_1$, $r_2$ and $h$. Let $S^*$ be the surface obtained from $S$ by
removing the images of the vertices of the rectangles. Identifying
$R$ with $R_1$ we can, proceeding as in the Abelian case,
write generators for the orbifold
fundamental group of $S^*$ as words in $h$, $r_1$ and $r_2$ or better
as words in $h$, $v=r_1r_2$ and $r_1$. We choose these words 
$w_1(h,v,r_1),\dots,w_s(h,v,r_1)$ to correspond to side pairings of 
the arrangement defining $S$.

Let $n$ be an integer satisfying the angle condition. 
Choose $A$ hyperbolic, $e_1$ and $e_2$ elliptic 
of order 2, as in section {\bf 1}, such that the group
$G_0=\langle A,e_1,e_2\rangle$ is discrete of signature $(0;2,2,2,n)$ or 
$(0;2,2,2,\infty)$ and such that  $\D/G_0-\{$point of order $n\}$ or 
$\D/G_0$, depending on the signature, is conformally equivalent to $S_0^*$. 

\proclaim{3.3 Proposition} Let $S$, the quadratic differential
$q$, the words $w_1,\dots,w_s$, $A$, $e_1$ and $e_2$ be as
above. 

Let $B=e_1e_2$ and let $G$ be the group generated by 
$w_1(A,B,e_1),\dots,w_s(A,B,e_1)$. Then properties {\rm (i)} to {\rm
(iv)} of {\bf 3.2} hold also in this case.
\endproclaim   

The proof follows exactly the same lines as the proof of {\bf 3.2}.
\bigskip

We end this section by describing the hyperbolic counterpart of the 
natural action of $\PSL_2(\Bbb Z)$ on the Teichm\"uller disk
orbit of a surface obtained from rectangles.

To avoid technical difficulties we need to restrict to surfaces for 
which the hyperbolic metric is non singular. To be more specific we
need
 
\proclaim{3.4 Definition} Let $(S,q)$ be a surface obtained from 
rectangles. We will say that $(S,q)$, or simply $S$, is balanced 
if all vertices correspond to zeros of the same order of the 
quadratic differential $q$.
\endproclaim

If $(S,q)$ is balanced then choosing the hyperbolic equiquadrangle
with interior angle $\pi/(n+2)$, where $n$ is the order of $q$ at the
vertices, leads to a non singular hyperbolic metric on the surface.

\proclaim{3.5} For the remainder of this section we will consider balanced 
surfaces with this hyperbolic metric in addition to the locally flat metric 
induced by the quadratic differential.
\endproclaim

Let $S$ be a surface obtained from rectangles, then $S$ has a natural 
decomposition into horizontal cylinders $\rs{C}_1,\dots,\rs{C}_p$ (see for 
example [Hu-Le]). If a cylinder $\rs{C}_i$ is formed of $n_i$ rectangles, 
we will say that it is of width $n_i$. By Lemma {\bf 2.1} or Lemma {\bf 2.2} 
the horizontal medians of these cylinders are disjoint simple closed 
geodesics $\gamma_i$. 

\proclaim{3.6 Proposition} Let $(S,q)$ be a balanced surface of genus $g$ 
obtained from squares, with the hyperbolic structure defined in {\bf 3.5}. 
Let $\rs{C}_1,\dots,\rs{C}_p$ be its decomposition into horizontal cylinders, 
$\rs{C}_i$ of width $n_i$ and let $\gamma_1,\dots,\gamma_p$ be the 
corresponding simple closed geodesics. Let 
$\gamma_{p+1}\dots,\gamma_{3g-3}$ be 
geodesics such that $\gamma_1,\dots,\gamma_{3g-3}$ defines a pants 
decomposition. Then there exist functions $f_{p+1},\dots,f_{3g-3}$,
defining lengths, and $tw_1,\dots,tw_{3g-3}$, defining twist parameters,
depending only on the combinatorics of the decomposition into squares 
and the choice of the $\gamma_j$, $j>p$, such that the surfaces in the 
$\PSL_2(\Bbb R)$ orbit of $(S,q)$ have Fenchel-Nielsen coordinates of the 
form
$$\multline \biggl(n_1\ell,tw_1(\ell)+\frac{t}{n_1},\dots,n_p\ell,tw_p(\ell)+
\frac{t}{n_p},\\ f_{p+1}(\ell),tw_{p+1}(\ell),\dots,f_{3g-3}(\ell),
tw_{3g-3}(\ell)\biggr)\endmultline\tag 3.6.1$$
with $\ell=2\,\ach(L)$, $L$ as in {\bf 2.1} and $t\in\Bbb R$.
\endproclaim

\noindent{\smc Proof.} We first consider the case when $(S,q)$ is obtained 
from rectangles. 

That the lengths of the $\gamma_i$, $1\leqslant i\leqslant k$,
is $n_i\ell$ follows immediately from the definition of the cylinders 
and the construction of the hyperbolic structure.

Let $\gamma_j$, $j>p$, be one of the additional geodesics. 
Because of our convention {\bf 3.5}, Proposition {\bf 3.2} (resp. {\bf 3.3}) 
implies that $\gamma_j$ is the image in $S$ of the axis of a word $w$ in
$A$ and $B$ (resp. $A$, $B$ and $e_1$) and the length of $\gamma_j$ is 
$2\ach(\trace(w)/2)$. 

Such a word depends only on the combinatorics of the decomposition into 
rectangles and the choice of $\gamma_j$. 
On the other hand the coefficients of $A$ and $B$ only depend on $n$, 
which is fixed, and $L=\cosh(\ell/2)$ (see Lemma {\bf 2.1}). In particular 
$\trace(w)$, and hence the hyperbolic 
length of $\gamma_j$, is a function of $L$.  
 
Let $\rs P_1$ and $\rs P_2$ be two pairs of pants pasted along 
$\gamma_i$. Then the twist parameter attached to $\gamma_i$ can be 
computed in terms of the lengths of boundary components of $\rs P_1$ and 
$\rs P_2$ and the length of a simple closed geodesic $\delta_i$ contained 
in $\rs P_1\cup\rs P_2$ (see for example [Bu], Proposition {\bf 3.3.12}).

Lift this geodesic to a curve $d_i$ in the arrangement of rectangles (which 
we assume as usual to have horizontal and vertical sides).
 
From the construction of this geodesic $\delta_i$ (see [Bu], Chap.~3, \S 3) 
the free homotopy class of $\delta_i$ in the surface only depends on the 
the combinatorics of the arrangement and the side identifications, and is
independent of the choice of the specific rectangle. Hence $\delta_i$ 
can be expressed as a fixed word in $A$, $B$ and $e_1$. The same argument
as above shows that the length of $\delta_i$ and the twist parameter for
$\gamma_i$ are functions of $L$.  
\smallskip
 
Let $\rs{R}_0$ be an equiquadrangle and let $A$ and $B$ be the
hyperbolic left-right and bottom-top side pairings. For each integer $j$ 
denote by $\rs R_j$ the equiquadrangle $A^j(\rs{R}_0)$. 

A cylinder $\rs C$ of width $n$ is the quotient of the
infinite union of the $\rs R_j$ by the identification of each $\rs R_j$
to $\rs R_{j+n+1}$ by $A^{n+1}$ (see left of figure 10 where we assume that 
the horizontal middle geodesic is the real axis in the unit disk). 

Now consider the
situation with the same $A$ but with $B$ replaced by $B_1=T\cdot B$, with 
$T$ as in Lemma {\bf 2.1} and with $Tw=\cosh(t\,\ach(L))$. There are two ways 
we can look at this situation, one is by replacing the equiquadrangles by 
the domain on the right of figure 7, another is to keep the 
lower half of the cylinder (in $\Im m(z)<0$) fixed and shift the upper 
half by $T$ (see right of figure 10). But since the length of the middle 
geodesic is $n\ell$ this amounts to applying a $\ach(Tw)/(n\ach(L))=t/n$ 
twist along the middle geodesic. This proves that for 
$1\leqslant j\leqslant k$ the twist parameter for $\gamma_j$ is now 
$tw_j(\ell)+t/n_j$.
\bigskip

\centerline{
 \SetLabels
\endSetLabels
 \AffixLabels{\psfig{figure=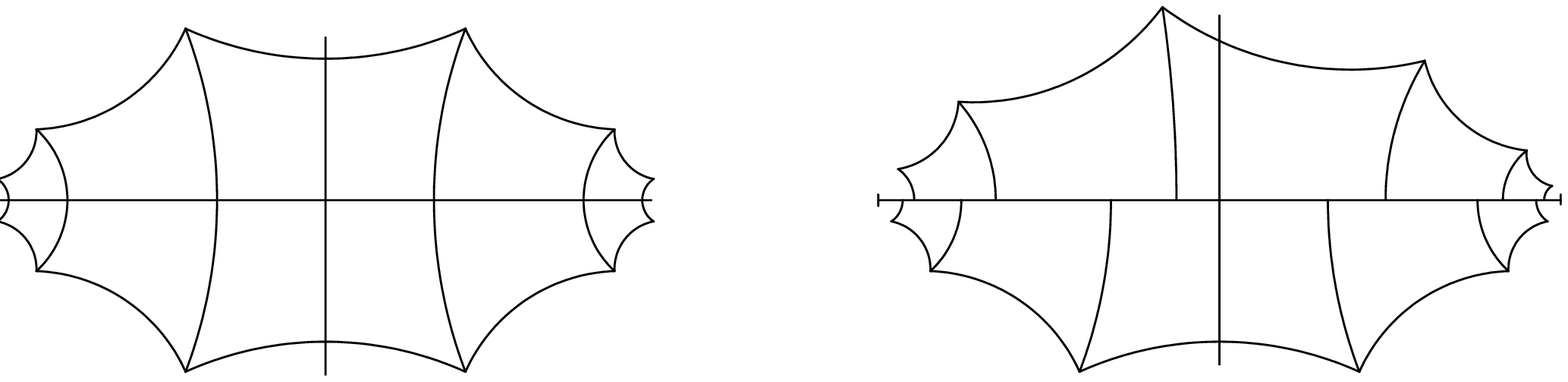,height=2.5cm}}}
\bigskip

\centerline{\bf Fig. 10}
\bigskip

To end the proof let $S$ be the surface obtained from the equiquadrangles
or otherwise said from $A$ and $B$ and let $S_1$ be the surface obtained 
from $A$ and $B_1=T\cdot B$. Let $S'$ and $S'_1$ be the surfaces with 
boundaries obtained from $S$ and $S_1$ respectively by cutting along 
$\gamma_1,\dots,\gamma_s$.
 The above construction shows that $S'$ and $S'_1$ are isometric. 
This shows that the remaining Fenchel-Nielsen coordinates, length of 
$\gamma_j$ and twist along $\gamma_j$, remain unchanged.
\bigskip
  
We started with a decomposition into horizontal cylinders, but we would 
have obtained the same type of result, starting with a decomposition into 
vertical cylinders. This would lead to another pants decomposition and 
another set of Fenchel-Nielsen coordinates
$$\multline \biggl(n'_1\ell',tw'_1(\ell')+\frac{t'}{n'_1},\dots,n'_{p'}\ell',
tw'_{p'}(\ell')+\frac{t'}{n'_{p'}},
\\ g_{p'+1}(\ell'),tw'_{p'+1}(\ell'),\dots,g_{3g-3}(\ell'),
tw'_{3g-3}(\ell')\biggr)\endmultline\tag 3.6.2$$ 
with $\ell'=2\,\ach(L')$, $L'$ as in {\bf 2.1}.

\proclaim{3.7 Corollary} Let $\rs D$ be a Teichm\"uller disk, 
$\PSL_2(\Bbb R)$ orbit of a surface obtained from squares.
Let 
$$\align &\left(n_1\ell,\,t_1,\dots,n_p\ell,\,t_p,
\ell_{p+1},\,t_{p+1},\dots,\ell_{3g-3},\,t_{3g-3}\right)\\ 
&\left(n'_1\ell',\,t'_1,\dots,n'_{p'}\ell',\,t'_{p'},
\ell'_{p'+1},\,t'_{p'+1},\dots,\ell'_{3g-3},\,t'_{3g-3}\right)\endalign$$
be pairs of Fenchel-Nielsen coordinates, of type $(3.6.1)$ and 
$(3.6.2)$, for surfaces in $\rs D$.

Let $\varphi_1$ be the composition of fractional Dehn-twists along 
$\gamma_1,\dots,\gamma_k$ that consists in replacing $t_1,\dots,t_p$ 
by $t_1+1/n_1,\dots,t_p+1/n_k$. Similarly let $\varphi_2$ 
be the composition of fractional Dehn-twists that consists in replacing 
$t'_1,\dots,t'_{p'}$ by $t'_1+1/n'_1,\dots,t'_{p'}+1/n'_{p'}$.

Then $\varphi_1$ and $\varphi_2$ generate the natural action of 
$\PSL_2(\Bbb Z)$ on $\rs D$.
\endproclaim

\noindent{\smc Proof.} Let $S$ be a surface in $\rs D$ tiled by
parallelograms with invariant $\tau$. The image of $S$ under 
$\left(\smallmatrix 1&1\\ 0&1\endsmallmatrix\right)$ is obtained by
replacing the initial parallelograms by those with invariant $1+\tau$. 
Let $A$  be the hyperbolic left-right side pairing and $B_1$ the other 
side pairing of the hyperbolic quadrangle associated to the
parallelogram defined by $\tau$. Replacing $\tau$ by $1+\tau$ is just
a change of basis of the lattice. By construction this corresponds to
replacing $A$ and $B_1$ by $A$ and $A\cdot B_1$. But this is just
taking $t=1$ in (3.6.1). This is the action of $\varphi_1$. 

The same argument shows that 
$\left(\smallmatrix 1&0\\ 1&1\endsmallmatrix\right)$ acts by
$\varphi_2$ on $\rs D$.
\bigskip

\noindent{\bf 3.8 Remark.} From the form of the matrices $A$ and $B$ in 
Lemma {\bf 2.1} and the proof of Proposition {\bf 3.6} one can easily 
deduce that all the $\cosh(f_j(\ell)/2)$ have algebraic expressions in $L$.

The same is true of the $\cosh(t_i\ell_i/2)$ since the exact expression 
given in [Bu] {\bf 3.3.11} only involves hyperbolic sines and hyperbolic 
cosines of lengths.
\bigskip

\subhead 4. The stairs and the escalator families \endsubhead
\medskip

In this section we will consider families constructed from the common 
admissible arrangement of rectangles illustrated in figure 11. Although 
these families belong to different strata we will see that the equations 
of the surfaces in these families are intimately related.
\medskip
\bigskip

\centerline{
 \SetLabels
\L(.04*.18) $a$\\
\L(.04*-.08) $b$\\
\L(.92*1.015) $c$\\
\L(.92*.765) $d$\\ 
\endSetLabels
 \AffixLabels{\psfig{figure=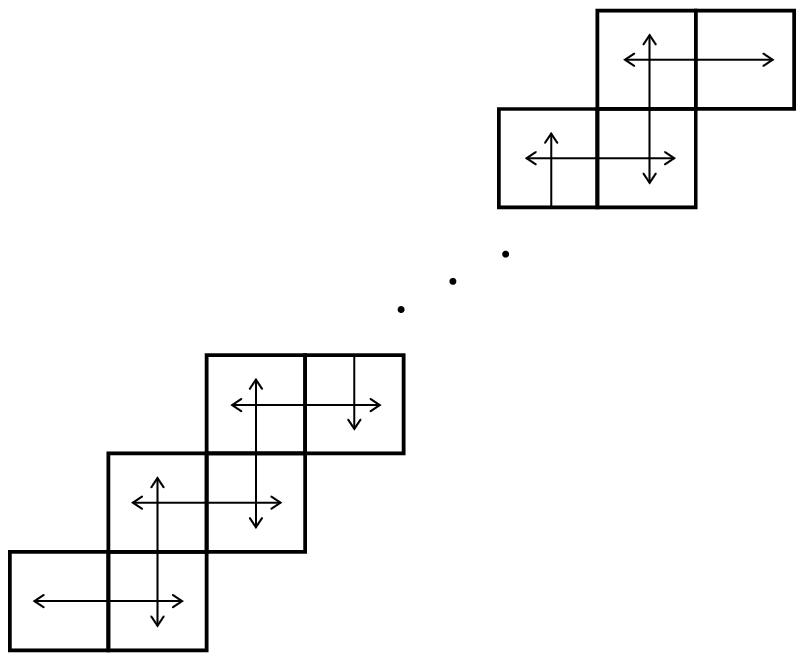,height=4cm}}}
\bigskip
\centerline{\bf Fig. 11}
\bigskip

The families are differentiated by the way the sides labeled $a$, $b$, $c$ 
and $d$ are identified. The families are as follows.
\medskip
 
St$_1(g)$: the number of rectangles is $2g$ and the identification 
is $a\sim b$, $c\sim d$. The surface is hyperelliptic of genus $g$ 
and the Weierstrass points are the centers of the rectangles  
and the midpoints of the edges labeled $a\sim b$ and $c\sim d$. The 
vertices map to two distinct points in the surface hence the surface 
is in the stratum $\Cal H(g-1,g-1)$. Surfaces in this 
family are often called stairs (see [Sch] or [M\"ol]).
\medskip

Esc$_1(g)$:  the number of rectangles is $2g+2$ with $g$ odd and the 
identification is $a\sim d$, $b\sim c$. The surface is hyperelliptic 
of genus $g$ and the Weierstrass points are the centers of the rectangles.
The vertices of the rectangles map to four distinct points hence the
surface is in the stratum $\Cal H((g-1)/2,(g-1)/2,(g-1)/2,(g-1)/2)$. 
By analogy with the stairs we call this family the escalator family. 
\medskip

Esc$_2$: the number of rectangles is $2g$ with $g$ odd and the 
identification is $a\sim d$, $b\sim c$. The surface is hyperelliptic 
of genus $g$ and the Weierstrass points are the centers of the rectangles 
and the vertices of the rectangles. The vertices of the rectangles map to 
two points and hence the surface is in the stratum $\Cal H(g-1,g-1)$. 
\medskip

Escb$_1(g)$:  the number of rectangles is $2g+2$ with $g$ even and the 
identification is $a\sim c$, $b\sim d$. The surface is hyperelliptic 
of genus $g$ and the Weierstrass points are the centers of the rectangles.
The vertices of the rectangles are mapped to four distinct points. Since 
the identifications $a\sim c$ and $b\sim d$ use a half-turn, the surface 
is in the stratum $\Cal Q((g-1),(g-1),(g-1),(g-1))$. The identification 
used here does not really deserve to be called of escalator type but we 
have named it by analogy with the Esc$_1$ family.
\medskip

Escb$_2(g)$: the number of rectangles is $2g$ with $g$ even and the 
identification is $a\sim c$, $b\sim d$. The surface is hyperelliptic 
of genus $g$ and the Weierstrass points are the centers of the 
rectangles and the vertices of the rectangles.
The vertices of the rectangles are mapped to two points and the surface 
is in the stratum $\Cal Q(2g-2,2g-2)$.
\medskip

Finally we also consider one last family based on a slightly different 
arrangement (see figure 12).
\bigskip

\centerline{
 \SetLabels
\endSetLabels
 \AffixLabels{\psfig{figure=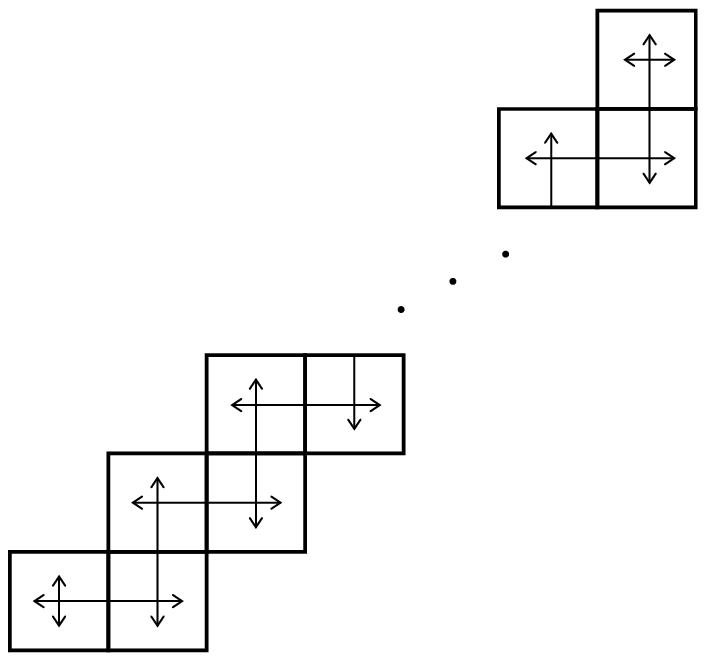,height=4cm}}}
\bigskip
\centerline{\bf Fig. 12}
\bigskip

St$_2(g)$ : the number of rectangles is $2g-1$. The surface is hyperelliptic 
of genus $g$ and the Weierstrass points are the vertices, all mapped to one 
point, the centers of the rectangles and the midpoints of the identified 
horizontal edges of the lower left rectangle and of the identified vertical 
edges of the upper right rectangle. The surface is in the stratum 
$\Cal H(2g-2)$.

Note that multiples of $2\,g$ satisfy the angle condition for the 
differential on surfaces in St$_1(g)$, Esc$_2(g)$ and Escb$_2(g)$, 
while multiples of $g+1$ satisfy the conditions for surfaces in Esc$_1(g)$ 
and Escb$_1(g)$. For surfaces in 
St$_2(g)$ we must choose multiples of $4\,g-2$. In particular we can 
replace the rectangles by equiquadrangles with angles $\pi/(2\,g)$ for
St$_1(g)$, Esc$_2(g)$ and Escb$_2(g)$, or angles $\pi/(g+1)$ for 
Esc$_1(g)$ and Escb$_1(g)$, or angles $\pi/(4\,g-2)$ for St$_2(g)$. If 
we do this we are under the hypothesis of {\bf 3.5} and the hyperbolic 
metric will be non-singular. To compute the corresponding Fuchsian groups 
we only need, by Propositions {\bf 3.2} and {\bf 3.3}, to compute  the words 
expressing the identifications in terms of the elementary horizontal and 
vertical identifications $h$ and $v$ of the rectangle and the rotation of 
angle $\pi$ at the center of the rectangle. This is taken care of by

 \proclaim{4.1 Lemma}
If the surface is in St$_1(g)$ the side pairings are given by,
$$\multline \langle\, v,(hv)^{g-1}hvh^{-1}(hv)^{1-g},(hv)^jh^2(hv)^{-j},
(hv)^ihv^2h^{-1}(hv)^{-i}\\ /\ 0\leqslant j\leqslant g-1,\
0\leqslant i\leqslant g-2\ \rangle
\endmultline\tag 4.1.1$$
If the surface is in St$_2(g)$ they are given by, 
$$\langle\ v,(hv)^{g-1}h(hv)^{1-g},(hv)^jh^2(hv)^{-j},
(hv)^jhv^2h^{-1}(hv)^{-j}\ /\ 
0\leqslant j\leqslant g-2\ \rangle\tag 4.1.2$$

For Esc$_1(g)$ and Esc$_2(g)$ the side parings are given by, 
$$\multline \langle\ (hv)^k,(hv)^{k-1}hv^{-1},(hv)^jh^2(hv)^{-j},
(hv)^ihv^2h^{-1}(hv)^{-i}\\ /\ 0\leqslant j\leqslant k-1,\
0\leqslant i\leqslant k-2\ \rangle\endmultline\tag 4.1.3$$
where $k=g+1$ for Esc$_1(g)$ and $k=g$ for Esc$_2(g)$. 

For Escb$_1(g)$ and Escb$_2(g)$ they are given by,
$$\multline\langle\  (hv)^kr, (hv)^{k-1}hr,(hv)^jh^2(hv)^{-j},
(hv)^ihv^2h^{-1}(hv)^{-i}\\ /\ 0\leqslant j\leqslant k-1,\
0\leqslant i\leqslant k-2\ \rangle\endmultline\tag 4.1.4$$
where again $k=g+1$ for Escb$_1(g)$ and $k=g$ for Escb$_2(g)$.
\endproclaim
 
For surfaces in St$_1$ and St$_2$ these computations have been done by 
G. Schmit\-h\"usen in [Sch]. Similar computations yield the result for 
the other families.
\bigskip

Our next objective is to compute equations for the associated algebraic 
curves. We do Esc$_1(g)$ and Escb$_1(g)$ first.

\proclaim{4.2 Proposition} Surfaces in Esc$_1(g)$ or Escb$_1(g)$ correspond 
to algebraic curves with an equation of the form
$$y^2=x^{2g+2}+a\,x^{g+1}+1\tag 4.2.1$$
(with $-2<a<2$ when the elementary tile is a rectangle). The Abelian
differential or the quadratic differential is a scalar multiple of
$$\omega=\frac{x^{(g-1)/2}dx}{y}\quad \text{ or }\quad 
q=\frac{x^{g-1}dx^2}{y^2}\ .\tag 4.2.2$$ 

Conversely if $C$ is a curve with equation $(4.2.1)$, with $a\neq\pm2$, then, 
if $g$ is odd, $(C,\omega)$ is in the $\SL_2(\Bbb R)$ orbit of surfaces
in Esc$_1(g)$ and if $g$ is even, $(C,q)$ is in the $\SL_2(\Bbb R)$ orbit of 
surfaces in Escb$_1(g)$.

Moreover the equation of the elliptic curve corresponding to one 
elementary rectangle (or more generally elementary parallelogram)
is
$$y^2=x^4+a\,x^2+1\tag 4.2.3$$ with the same $a$ as in $(4.2.1)$
and we have
$$a=\frac{2\,\mu-4}{\mu}\ ,$$
where $\mu$ is the invariant of the rectangle as defined in section {\bf 2}.
\endproclaim

\noindent{\smc Proof.} We first note that, 
from the combinatorial structure of the identifications and the fact 
that the angle at the vertices of the equiquadrangles, defining the 
non-singular hyperbolic metric, is $\pi/(g+1)$ one 
can easily check that the surface has an automorphism $f_g$ of order $g+1$, 
the fixed points of which are the four images of the vertices.

Now assume that $g$ is odd and hence the surface $S$ is in Esc$_1(g)$. The 
quotient of $S$ under ${f_g}^{2}$ is a surface $E$ in Esc$_1(1)$ 
obtained from four rectangles. The quotient map is ramified precisely at the 
vertices of the rectangles, i.e.~ the fixed points of $f_g$.

To reconstruct this situation, we first note that in this case $f_1$ is an 
involution and we can always assume that it is induced by $x\mapsto -x$.
The four vertices will then be points above $x=0$ and $x=\infty$.
Label the rectangles from 1 to 4 starting with the lower left. The choice 
we have just made implies that if the $x$ coordinates of the center of
rectangles 1 and 2 are $x_1$ and $x_2$ then the centers of rectangles 3 
and 4 are $-x_1$ and $-x_2$. We can still make one choice,so we choose the
midpoints of the horizontal edges of rectangles 2 and 3 to have $x$-coordinate 
1. The midpoints of the horizontal edges of rectangles 1 and 4 will then 
correspond to $-1$. But this implies that the involution obtained by rotating 
the arrangement of rectangles by angle $\pi$ will be induced by $x\mapsto 1/x$.
In which case we have $x_2=-1/x_1$. Hence an equation of the form
$$y^2=(x^2-x_1^2)(x^2-1/x_1^2)=x^4+a\,x^2+1\tag 4.3$$
Now we have four involutions induced by half-turn around:
\roster
\item the vertices; this is
$f_1$ which is $(x,y)\mapsto (-x,y)$;
\item the midpoints of
horizontal sides; this is $(x,y)\mapsto (1/x,y/x^2)$;
\item the centers of the rectangles; this is
$(x,y)\mapsto (x,-y)$;
\item the midpoints of vertical edges.
\endroster
Since (4) is the composition of (1), (2) and (3), it is $(x,y)\mapsto
(-1/x,-y/x^2)$, and hence its fixed points are the points with
$x$-coordinate $\pm i$.

To end the description of $E$ we note that for rectangles we have real 
structures obtained by taking reflections in the horizontal or vertical 
medians of the rectangles. The first must fix the points above $\pm i$ and
the second the points above $\pm 1$. This implies that the first is 
$(x,y)\mapsto (1/\bar{x},-\bar{y}/\bar{x}^2)$ and the second is 
$(x,y)\mapsto (1/\bar{x},\bar{y}/\bar{x}^2)$. This in turn implies that 
$|x_1|=1$ and hence $-2<a<2$.

Now consider the genus $g$ curve with equation 
$$y^2=x^{2g+2}+a\,x^{g+1}+1\ ,\tag 4.4$$
$g=2\,n+1$ odd. This curve has an obvious automorphism of order $g+1$ defined
by $\varphi:(x,y)\mapsto (\zeta\,x,y)$, where $\zeta$ is a primitive 
$(g+1)$-th root of unity. The fixed points of $\varphi$ are the points above 
$x=0$ and the points at infinity. The quotient of this curve under  
$\varphi^2$ 
is the curve with equation (4.3). The quotient 
morphism is $(x,y)\mapsto (x^{n+1},y)$ which is precisely ramified at the 
points above $x=0$ and the points at infinity. This proves (4.2.1) for 
surfaces in Esc$_1(g)$.
\smallskip

If $g$ is even and $S$ is in Escb$_1(g)$ we consider the surface $S'$ in
Esc$_1(2\,g+1)$ obtained from the same rectangles. Then $S'$ is a double 
cover of $S$ and the covering is precisely ramified at the vertices of the 
rectangles. But we know from the above that $S'$ has an equation of the 
form $y^2=x^{4\,g+4}+a\,x^{2g+2}+1$. From this it is not hard to deduce 
that $S$ has again an equation of the form (4.2.1).
\smallskip

Since the differentials vanish at the points above $x=0$ and the points at
infinity, we also get (4.2.2).
\smallskip

The third assertion follows from the fact that the image in moduli space 
of the set of surfaces in Esc$_1(g)$ (resp. Escb$_1(g)$) is an 
algebraic curve.

To end the proof we only need to note that the surface $E$ obtained from 
four rectangles is obviously isomorphic to the one obtained from one 
rectangle and that the curve defined by (4.3) is isomorphic to
$$y^2=x\,(x-1)\,\left(x-\frac4{2-a}\right)$$
the isomorphism being induced by 
$$x\mapsto \frac{2\,x_1(x-x_1)}{(1+x_1^2)(x_1x-1)}\ .$$ 
\medskip

If the rectangles are in fact squares we have $\mu=2$ hence,

\proclaim{4.5 Corollary} Let $S$ be a surface of genus $g$ in Esc$_1$ or 
Escb$_1$. If $S$ is tiled by squares, then an equation for $S$ is 
$$y^2=x^{2g+2}+1\ .$$
\endproclaim

Form {\bf 4.2} we are going to deduce the other cases. We do  
Esc$_2(g)$ and Escb$_2(g)$ first.

\proclaim{4.6 Proposition} Let $S$ be a surface of odd, respectively even,
genus $g$ in Esc$_2$, respectively Escb$_2$. Then the corresponding algebraic
curve has an equation of the form
$$y^2=x\,(x^{2g}+a\,x^g+1)\ .\tag 4.6.1$$
The Abelian, resp.\ quadratic, differential is a scalar multiple of
$$\omega=\frac{x^{(g-1)/2}dx}{y}\quad \text{ resp.\ }\quad
q=\frac{x^{g-1}dx^2}{y^2}\ .\tag 4.6.2$$

Conversely if $C$ is a curve with equation $(4.6.1)$, with
$a\neq\pm2$, 
then, if $g$ is odd, $(C,\omega)$ is in the $\SL_2(\Bbb R)$ orbit of surfaces
in Esc$_2$ and if $g$ is even, $(C,q)$ is in the $\SL_2(\Bbb R)$ orbit of 
surfaces in Escb$_2$.

Moreover the equation of the elliptic curve corresponding to one 
elementary rectangle (or more generally elementary parallelogram)
is
$$y^2=x^4+a\,x^2+1\tag 4.6.3$$
with the same $a$ as in $(4.6.1)$ and this parameter only depends on the 
elementary parallelogram used in the construction.

In particular if $S$ is tiled by squares, then an equation for $S$ is
$$y^2=x\,(x^{2g}+1)\ .\tag 4.6.4$$
\endproclaim

\noindent{\smc Proof.} Let $S$ be in Esc$_2(g)$ or in Escb$_2(g)$ depending
on the parity of $g$. Consider the surface $S'$ in Esc$_1(2\,g-1)$ obtained 
from the same rectangles. Then $S'$ is an unramified double cover of $S$. 
More precisely let $h$ be $f_{2g-1}^{g}$ composed with the hyperelliptic
involution. Then $S$ is the quotient $S'/h$. From {\bf 4.2} we know that $S'$ 
has an equation of the form (4.2.1) and that $h$ is $(x,y)\mapsto (-x,-y)$. 

From this it is immediate to deduce (4.6.1) the quotient map being 
$(x,y)\mapsto (x^2,x\,y)$. The rest easily follows from 
Proposition {\bf 4.2} and its proof.
\bigskip

For St$_1(g)$ and St$_2(g)$ the expressions we find are not so nice, but
they can nevertheless be deduced from {\bf 4.2} and {\bf 4.6}. 

Let $S'$ be a surface in Esc$_1(2\,g-1)$. The involution induced by rotation 
of the arrangement of rectangles has 4 fixed points, but if we compose this 
involution with the hyperelliptic involution we obtain an involution with 
no fixed points. The quotient $S$ of $S'$ under this last involution is in 
St$_1(g)$ and obtained from the same rectangles. We note also that all 
surfaces in St$_1(g)$ are obtained in this way.

On the other hand, by Proposition {\bf 4.2}, $S'$ has an equation of the form 
$y^2=x^{4g}+a\,x^{2g}+1$ and the involution we are considering is 
$(x,y)\mapsto(1/x,y/x^{2g})$. An equation for the quotient is not as 
easy to express as in previous cases but we can do the following. Let
$x_1,1/x_1,\dots,x_{2g},1/x_{2g}$ be the roots of $x^{4g}+a\,x^{2g}+1$.
Let $t_k=i\frac{1+x_k}{1-x_k}$. Then $y^2=\prod(x^2-t_k^2)$ is also an 
equation for $S$ and the involution $(x,y)\mapsto (1/x,y/x^{2g})$ is now 
$(x,y)\mapsto (-x,-y)$. Hence

\proclaim{4.7 Proposition} Let $S$ be in St$_1(g)$ and let $t_i$ be as 
above. Then an equation for $S$ is
$$y^2=x\,\prod(x-t_k^2)\ .$$
The differential defining the locally flat metric on $S$ is 
$$\frac{(x+1)^{g-1}dx}y\ .$$
\endproclaim 

For surfaces in St$_2(g)$ we can use the same argument but starting with 
a surface $S'$ in Esc$_2(2\,g-1)$. Let 
$x_1,1/x_1,\dots,x_{2g-1},1/x_{2g-1}$ be the roots of 
$x^{4g-2}+a\,x^{2g-1}+1$ and let $t_k=i\frac{1+x_k}{1-x_k}$.

\proclaim{4.8 Proposition} Let $S$ be in St$_2(g)$ and let the $t_k$ be
as above, then an equation for $S$ is
$$y^2=x(x+1)\prod(x-t_k^2)\ .$$
The differential defining the locally flat metric on $S$ is
$$\frac{(x+1)^{g-1}dx}y\ .$$
\endproclaim 

If the surfaces are tiled by squares, i.e.~ if we have $a=0$ in the 
equations for $S'$ in Esc$_1(2g-1)$ or in Esc$_2(2g-1)$, a tedious but 
elementary computation shows that we have

\proclaim{4.9 Corollary} If $S$ is the surface in St$_1(g)$ tiled 
by squares then $S$ has for equation
$$y^2=x\left(\sum_k (-1)^k\binom{4g}{2k} x^{k}\right)\ .$$ 

If $S$ is the surface in St$_2(g)$ tiled 
by squares then $S$ has for equation
 $$y^2=x(x+1)\left(\sum_k (-1)^{k+1}\binom{4g-2}{2k} x^{k}\right)
\ .$$
\endproclaim

\subhead 5. Balanced genus 2 surfaces tiled by four rectangles \endsubhead
\medskip

In the last section we have seen relations between different 
surfaces tiled by rectangles. The aim of this section is to explore in
more detail the consequences of Proposition {\bf 3.6} and Corollary
{\bf 3.7} in the case of balanced surfaces of genus 2 tiled by four 
rectangles. We will also exhibit the action of other fractional Dehn twists 
exchanging the different families.

\proclaim{5.1 Proposition} There are exactly four $\PSL_2(\R)$ orbits 
of balanced surfaces tiled by four rectangles. These are the orbits 
of the surfaces described in figure $13$ {\rm (}where identifications
are indicated by numbers and top-top or bottom-bottom
identifications are by half-turns while all others are by translations{\rm)}.
\endproclaim

\noindent{\smc Proof.} Consider a balanced surface tiled by 4 rectangles.  
We leave aside the
case when the surface is a torus.  The angles at the vertices of the
rectangles add up to $4 \cdot 2 \, \pi = 8 \, \pi$.  Thus the vertices
cannot be identified to a single point on the surface: this would be a
zero of order $6$ for the quadratic differential, but the
multiplicities of the zeros of a quadratic differential add up to
$4\,g - 4$, a multiple of $4$.  The vertices can be identified to two
points on the surface, each of angle $4\,\pi$, corresponding to two
zeros of order $2$ for the quadratic differential (possibly the square
of an Abelian differential with two zeros of order $1$).  The case
when the vertices are identified to four points on the surface has
been ruled out (they would be points of angle $2\,\pi$ and the surface
would be a torus).

Let us therefore enumerate all (connected) surfaces obtained from $4$
rectangles, the vertices of the rectangles being identified to $2$
points of angle $4\,\pi$.  Decomposing the surface into horizontal 
cylinders, we find one of the four possible situations: all $4$ 
rectangles are lined up horizontally and form a single horizontal 
cylinder of width $4$ rectangles; they form $2$ cylinders of widths 
$2$ and $2$ or of widths $1$ and $3$; or they form three cylinders of 
widths $1$, $1$, and $2$.  There cannot be $4$ horizontal cylinders 
of widths $1$: if all rectangles have their left and right sides 
glued together then the surface is a torus.

Now consider the four possible cylinder decompositions in turn.

\noindent{\bf Case $1$}: one cylinder of width $4$.  Consider the $4$
rectangles lined up horizontally, forming a wide rectangle; the
surface results from the identification of the left and right sides of
this wide rectangle, and of pairwise identifications of the horizontal
rectangle edges.  Label the bottom sides of the rectangles 1, 2,
3, 4, from left to right, and likewise their top sides 5, 6,
7, 8.

Gluings 1--2, 1--4, 2--3, 3--4, 5--6, 5--8, 
6--7, 7--8 are forbidden, they would yield a cone point of 
angle $\pi$.  Gluing 1--5 and 2--6 simultaneously is also 
forbidden, it would yield a point of angle $2\,\pi$.  Avoiding such 
gluings, let us enumerate the possible surfaces we can obtain by 
gluing pairs of sides.

(1--3, 2--4, 5--7, 6--8), (1--3, 2--5, 4--7, 6--8), (1--3, 2--6, 4--8, 5--7),
(1--3, 2--7, 4--5, 6--8) and (1--3, 2--8, 4--6, 5--7) give five different
surfaces in $\Cal Q (2,2)$.

(1--5, 2--8, 3--7, 4--8) and (1--6, 2--5, 3--8, 4--7) give two different
surfaces in $\Cal H (1,1)$.

The remaining possible gluings, starting with 1--7 or 1--8, would yield
surfaces already listed (one can see that by cutting the leftmost
square and pasting it to the right, or cutting the rightmost square
and pasting it to the left).

\noindent{\bf Case 2}: two cylinders of widths $2$ and $2$.

The two cylinders have to be glued one to the other, so let us number
the rectangles $R_1$, $R_2$, $R_3$, $R_4$, and
suppose we start with $R_2$ glued to the right of $R_1$,
$R_3$ glued on top of $R_2$, $R_4$ glued to the right
of $R_3$.  Label 1, 2, 4 the bottom sides of $R_1$,
$R_2$, $R_4$, and 3, 5, 6 the top sides of $R_1$,
$R_3$, $R_4$.

Gluings 1--2, 5--6 are forbidden, they would yield a cone point of angle
$\pi$.  Gluing 1--5 and 2--6 simultaneously is also forbidden, it would
yield a point of angle $2\,\pi$.

This leaves the following possible gluings.

(1--3, 2--5, 4--6) and (1--3, 2--6, 4--5) give surfaces in $\Cal H (1,1)$,
(1--4, 2--5, 3--6), (1--4, 2--6, 3--5) give surfaces in $\Cal Q (2,2)$.
(1--5, 2--3, 4--6) and (1--5, 2--4, 3--6) give the surfaces already listed
as (1--3, 2--6, 4--5) and (1--4, 2--6, 3--5) respectively. (1--6, 2--3, 4--5)
gives a surface in $\Cal H (1,1)$, and (1--6, 2--4, 3--5) a surface in 
$\Cal Q (2,2)$.

\noindent{\bf Case 3}: two cylinders of widths $3$ and $1$.

Again let us number the rectangles $R_1$, $R_2$, $R_3$,
$R_4$, and without loss of generality start with $R_2$ glued
to the right of $R_1$, $R_3$ glued to the right of $R_2$, $R_4$ 
glued on top of $R_1$.  Label 1, 2, 3 the bottom sides of $R_1$, 
$R_2$, $R_3$, and 4, 5, 6 the top sides of $R_4$, $R_2$, $R_3$.

The following gluings are possible, and give rise to different 
surfaces.

(1--4, 2--6, 3--5), (1--5, 2--6, 3--4), (1--6, 2--5, 3--4).

\noindent{\bf Case 4}: three cylinders of widths $1$, $1$ and $2$.

Without loss of generality we can start with $R_2$ glued to the
right of $R_1$, $R_3$ glued below $R_1$, and $R_4$
glued either on top of $R_1$ (sub-case 1) or on top of $R_2$
(sub-case 2).

In sub-case 1, label 1 and 2 the bottom sides of $R_4$ and $R_2$, 
3 and 4 the top sides of $R_3$ and $R_2$.  The
possible gluings are (1--2, 3--4) and (1--4, 2--3).

In sub-case 2, label 1 and 2 the bottom sides of $R_4$ and $R_2$, 
3 and 4 the top sides of $R_1$ and $R_3$.  The
possible gluings are (1--2, 3--4) and (1--3, 2--4).

The gluings (1--3, 2--4) of sub-cases 1 and 2 yield the same surface,
so we get three different surfaces from case 4.

\medskip

Our case study evidenced 19 different balanced surfaces tiled with 4
rectangles; one can compute their $\SL(2,{\Z})$-orbits and see that
these 19 surfaces fall into 4 different orbits, representatives of
which are presented in figure 13 (the detailed description of the 
19 different cases is the object of the rest of this section).
\bigskip

\centerline{
 \SetLabels
\L(.24*.54) $(A)$\\
\L(.74*.54) $(B)$\\
\L(.385*0) $(C)$\\
\L(.885*0) $(D)$\\
\endSetLabels
 \AffixLabels{\psfig{figure=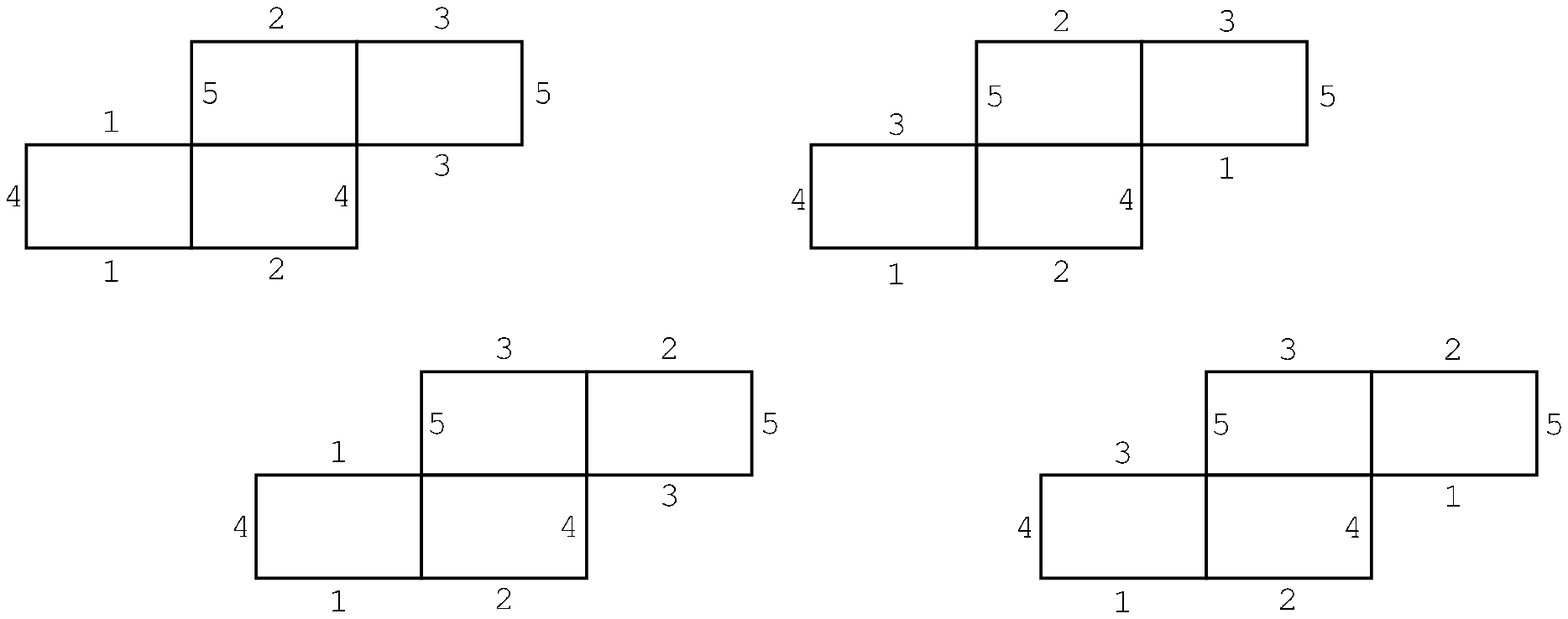,height=5cm}}}
\bigskip
\centerline{\bf Fig. 13}
\bigskip

To give a full description of these families we need a more detailed 
description of the geometry of the hyperbolic quadrangles we will 
consider.

Let $L$ be as in Lemma {\bf 2.1}, then since the angle is $\pi/4$
we have $L'=\sqrt{\dfrac{2\,L^2-1}{2\,L^2-2}}$. If we take $A$ and $B$ 
as in (2.1.1), with this value of $L'$, and $T$ as in (2.1.2),
then $A$ and $T\cdot B$ generate a Fuchsian group of signature
$(1;2)$. Moreover all Fuchsian groups of signature $(1;2)$ are
conjugate to one of this form.
Note also that in this context the elliptic element $e_1$ of order two 
at the center of symmetry of the quadrangle is the conjugate of  
$z\mapsto -z$ by a matrix of the same form as $T$ but with $Tw$ replaced by 
$\sqrt{(Tw+1)/2}$.
 
Call $L_1$ (resp. $L_2$) the hyperbolic cosine of the hyperbolic distance 
between $p_1$ and $p_3$ (resp. $p_2$ and $p_4$) in figure 7 (left), call 
$L_3$ (resp. $L_4$)  the hyperbolic cosine of the hyperbolic distance 
between $q_1$ and $q_2$ (resp. $q_1$ and $q_4$) (figure 7 left) and 
finally call $L'_2$ (resp. $L'_4$)  the hyperbolic cosine of the hyperbolic 
distance between $p_3$ and $p'_2$ (resp. $q_4$ and $q'_1$) (figure 7 right).

\proclaim{5.2 Lemma} We have the following relations
\roster
\item"(i)" $L_2=\dfrac{L_1+1}{L_1-1}$;
\item"(ii)" $L_3=2\,L_1+1$;
\item"(iii)" $L_4=2\,L_2+1$;
\item"(iv)" $L'_2=4\,Tw^2(L_2+1)-1=Tw^2\dfrac{2\,L_1}{L_1-1}-1$ $($with 
$Tw$ as in $(2.1.2)$$)$;
\item"(v)" $L'_4=2\,L'_2+1$.
\endroster
\endproclaim

\noindent{\smc Proof.} The first three assertions are immediate consequence of 
the relations between the side lengths of a trirectangular quadrangle with
remaining angle $\pi/4$ (see [Bu], p.454). For (iv) we note that the trace of
$T\cdot B$ is $L''=2\,Tw\,L'$ ($L'$ as above). But $L_2=2\,{L'}^2-1$ and 
$L'_2=2\,{L''}^2-1$, combined with (i) this yields (iv). For (v) we note that
$q'_1$ is the image of $q_4$ under the action $T\cdot B$ and from this the 
distance between $q_4$ and $q'_1$ can easily be computed and the relation 
checked.
   
\proclaim{5.3 Proposition} Let $S_A$, $S_B$, $S_C$ or $S_D$ be in the 
$\SL_2(\R)$ orbit of one of the surfaces of type $(A)$, $(B)$, $(C)$ or 
$(D)$ of {\bf 5.1}. Then Fenchel-Nielsen coordinates of these surfaces are 
of the form
\roster
\item"(i)" $(\ell,tw,\ell,tw,\ell',0)$ for $S_A$;
\item"(ii)" $(\ell,tw,\ell,tw,\ell',\frac12)$ for $S_B$;
\item"(iii)" $(\ell,tw+\frac12,\ell,tw,\ell',0)$ for $S_C$;
\item"(iv)" $(\ell,tw+\frac12,\ell,tw,\ell',\frac12)$ for $S_D$,
\endroster
where $\cosh(\ell'/2)=2\,\cosh(\ell/2)+1$.
\endproclaim

\noindent{\smc Proof.} We first note that the surfaces are balanced and 
that in all cases the quadratic differential defining the locally flat 
metric has two zeros of order 2. Hence replacing the rectangles by an  
equiquadrangle with interior angle $\pi/4$ (left of figure 7) leads to a 
smooth hyperbolic surface. Moreover since the angles at the vertices are 
$\pi/4$ the union of the arcs labeled 1 and 3 in figure 13 (all cases) 
forms a simple closed hyperbolic geodesic that we will call $\gamma_3$. 

The medians of the horizontal cylinders are also simple closed hyperbolic 
geod\-esics $\gamma_1$ (for the lower cylinder) and $\gamma_2$ (for the upper 
cylinder). Since the $\gamma_i$ do not intersect they define pants 
decomposition of the surfaces.  

To obtain Fenchel-Nielsen coordinates we start with case $(A)$ for 
equiquadrangles. In this case reflection along $\gamma_1$ clearly fixes 
point-wise $\gamma_2$ and $\gamma_3$, from this it follows that the twist 
parameters are all zero in this case hence Fenchel-Nielsen coordinates of
the form   $(\ell,0,\ell,0,\ell',0)$. The assertion that   
$\cosh(\ell'/2)=2\,\cosh(\ell/2)+1$ immediately follows from {\bf 5.2} (ii).

We have here only considered the euquiquadrangle case to use the fact that
the $\gamma_i$ are multigeodesic, but clearly, using {\bf 5.2} (iv), we can 
replace the rectangles by more general quadrangles of the form illustrated 
on the right of figure 7, this yields Fenchel-Nielsen coordinates of the
form  $(\ell,tw,\ell,tw,\ell',0)$ with of course again 
$\cosh(\ell'/2)=2\,\cosh(\ell/2)+1$. That this describes the full 
$\SL_2(\R)$ orbit follows from Proposition {\bf 3.6} (for this case
see also [Si2] section {\bf 3}).

The claim for surface of type $(B)$ can also be deduced form [Si2] 
section {\bf 3}, but to cover also the relation between type $(C)$ and
$(D)$ we are going to deduce this from a more general result.

\proclaim{5.4 Lemma} Let $S$ be a surface in Esc$_2(g)$, 
if $g$ is odd, or Escb$_2(g)$, if $g$ is even. If $S$ is in Esc$_2(g)$
then the union of the arcs labeled $a\sim d$ and $b\sim c$ (see figure
11) forms a simple closed hyperbolic geodesic $\gamma$. If $S$ is in  
Escb$_2(g)$ then the union of the arcs labeled $a\sim c$ and $b\sim d$ 
forms a simple closed hyperbolic geodesic that we will also denote by 
$\gamma$.

Let $S'$ in St$_1(g)$ be obtained from the same rectangles.  In the same way  
$\{a\sim b\}\cup\{c\sim d\}$ defines a simple closed geodesic $\gamma'$ 
in $S'$. 

Then $S'$ is obtained from $S$ by applying a half-Dehn twist along the 
geodesic $\gamma$. Conversely $S$ is obtained from $S'$ by applying a 
half-Dehn twist along the geodesic $\gamma'$.
\endproclaim

\noindent{\smc Proof.}  Label consecutively the rectangles of figure 11
from 1 to $2g$ starting with the rectangle on the lower left hand
side. Note also that both for $S$ and $S'$  the vertices lie in two orbits. 

In all cases label $\circ$ those in the 
orbit of the upper left corner of rectangle 1 and label 
$\bullet$ those in the other orbit. 

Since the angle at the vertices is $\pi/(2g)$
the combinatorics of the rectangles at these points is described in figure 
14. In the upper figure the geodesic $\gamma$ corresponds to the horizontal 
line passing through $\circ$ and $\bullet$, while in the lower figure the 
horizontal line represents $\gamma'$. This proves the assertions on
$\gamma$  and $\gamma'$.

Since $\circ$ and $\bullet$ split $\gamma$ into two arcs of equal 
hyperbolic lengths it follows from the defintion of a half-Dehn 
twist that we pass from the upper part of figure 14 to the lower part 
by performing a half-Dehn twist along $\gamma$. But the upper part gives the 
combinatorics for Esc$_2$ or Escb$_2$ while the lower gives the combinatorics 
for St$_1$. This proves Lemma {\bf 5.4}.
\bigskip

\centerline{
 \SetLabels
\L(.39*.5) Esc$_2$ or Escb$_2$\\
\L(.49*-.08) St$_1$\\
\endSetLabels
 \AffixLabels{\psfig{figure=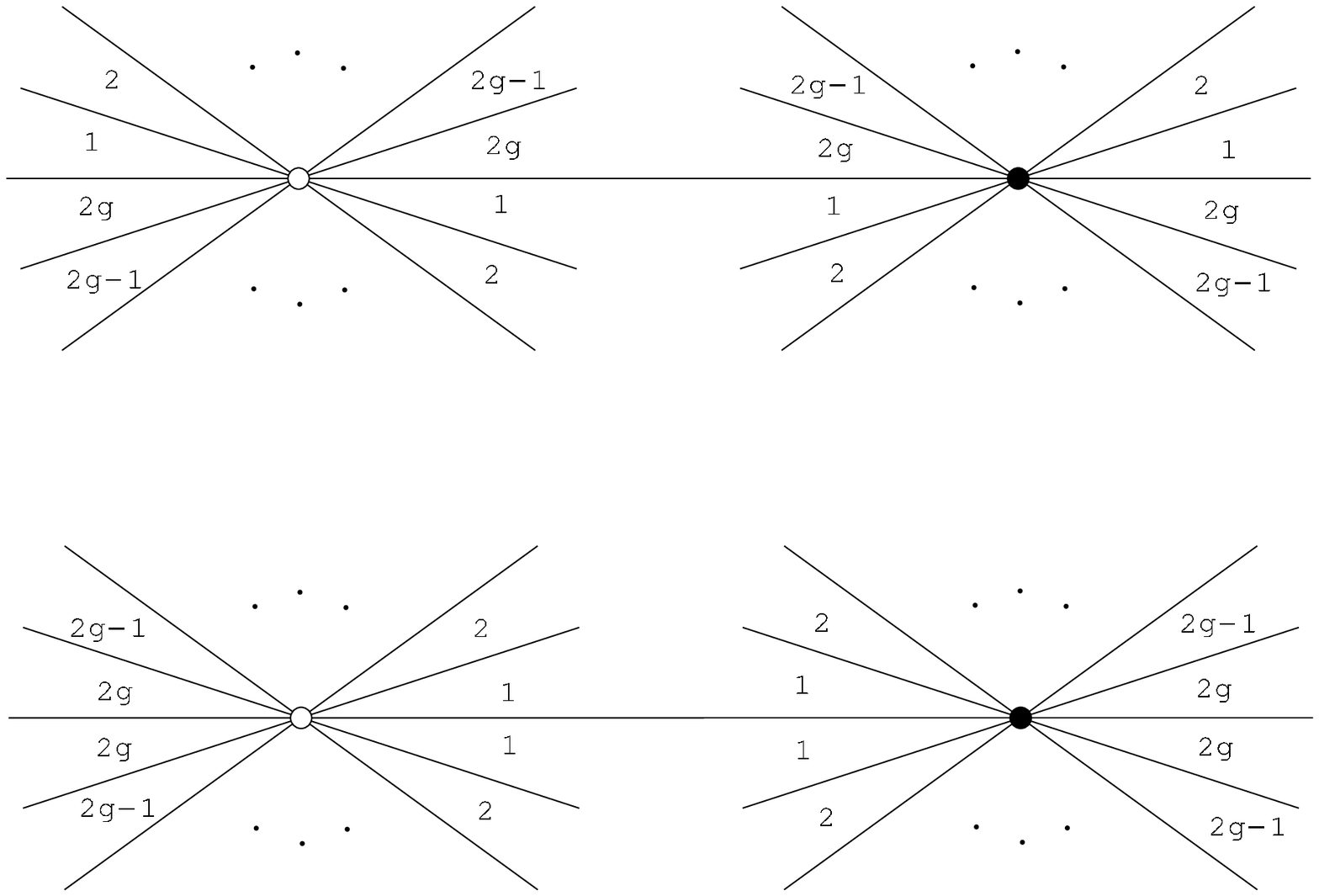,height=6cm}}}
\bigskip
\bigskip

\centerline{\bf Fig. 14}
\bigskip  

For simplicity we have formulated {\bf 5.4} for rectangles and 
equiquadrangles, but obviously, althoug the geodesics may not conicide, 
the statement can be generalized for parallelograms by Poposition 
{\bf 3.6}.

Since $S_A$ is in St$_1(2)$ and $S_B$ in Escb$_1(2)$ the statement of 
{\bf 5.3} for type $(B)$ follows. 

The same argument also shows that one passes from $S_C$ to $S_D$ by a 
half-Dehn twist along $\gamma_3$.

To end the proof we note that $S_C$ (resp. $S_D$) is obtained from $S_A$ 
(resp. $S_B$) by applying a half-Dehn twist along $\gamma_2$, as can
be immediately checked by looking at the identifications.
\bigskip

We end this section by computing the equations for the different
families.
\medskip

\noindent{\bf Case $(A)$}. This family is of course St$_1(2)$ and equations
can be recovered by applying the results of section {\bf 2}. But to 
highlight the links between the different cases we are going to use
a slightly different approach.

We use the fact that the surface has a non hyperelliptic involution 
induced by a rotation of angle $\pi$. Hence we can look for an equation
of the form 
$$y^2=(x^2-a^2)(x^2-1)(x^2-b^2)\ .\tag 5.5$$
If we label the rectangles $R_1$ to $R_4$ starting with the lower left one, 
the existence of the additional involution imposes that the 
midpoints of the horizontal edges of rectangles $R_2$ and $R_3$ have 
$x$-coordinate $0$ or $\infty$. We choose 0.
Since composing this involution with the hyperelliptic one fixes the vertices, 
these will have $x$-coordinate $\infty$. We normalize further by 
chosing the $x$-coordinate of the Weierstrass point at the center of rectangle 
$R_1$ to be $-1$ which forces the center of the rectangle $R_4$ to have 
$x$-coordinate $1$. Finally we choose the Weierstrass points $(-b,0)$ and
$(b,0)$ to be the midpoints of the horizontal edges of rectangles
$R_1$ and $R_4$ 
respectively. With this fixed, a map from the surface to the genus 1 surface 
tiled by one rectangle is induced by the map 
$$f:x\mapsto \frac{x^2(x^2-b^2)}{1-b^2}\ .\tag 5.6$$
But now $a$ and $-a$ are simply solutions of $f(x)=1$ distinct from $\pm1$. 
Hence 
$$a=\sqrt{b^2-1}\ .\tag 5.7$$ 
To complete the description of the genus 1 quotient 
(or alternatively recover $b$ in terms of the genus 1 quotient) we 
only need to compute for which values of $\lambda$ the equation
$x^2(x^2-b^2)-\lambda\,(1-b^2)=0$ has a double root. This yields
$$\lambda=\frac{b^4}{4\,(b^2-1)}\ .\tag 5.8$$  
This does not conform to our convention on the $\mu$ invariant for genus 1
but we easily find,
$$\mu=\frac{\lambda}{\lambda-1}=\frac{b^4}{(b^2-2)^2}=
\frac{(a^2+1)^2}{(a^2-1)^2}\ .\tag 5.9$$

Summarizing, we have an equation for the algebraic curve, of the form 
$$y^2=(x^2-a^2)(x^2-1)(x^2-a^2-1)\tag 5.10$$
and a degree 4 map, ramified at the points at infinity,
$$(x,y)\mapsto\left(\frac{(2\,x^2-a^2-1)^2}{(a^2-1)^2},
y\frac{4\,x(2\,x^2-a^2-1)}{(a^2-1)^3}\right)\tag 5.11 $$
onto the genus 1 curve defined by
$$y^2=x(x-1)\left(x-\frac{(a^2+1)^2}{(a^2-1)^2}\right)\ .\tag 5.12$$

To obtain the values of $a$ for the other combinations of rectangles in 
the same $\SL_2(\Bbb Z)$-orbit we note that the difference between
cases $(A_1)$ to $(A_6)$ (see figure 15)
is in the repartition of the Weierstrass points among vertices, 
centers of rectangles, horizontal edges and vertical edges. These are 
$(0,4,2,0)$ $(A_1)$, $(0,0,2,4)$ $(A_2)$, $(0,2,4,0)$ $(A_3)$, $(0,4,0,2)$ 
$(A_4)$, $(0,0,4,2)$ $(A_5)$ or $(0,2,0,4)$ $(A_6)$ (see figure 15).
We can then procceed as above.

But a better solution here is to replace, in equation (5.9), $\mu$ by
$1-\mu$, $1/\mu$ and so forth. For this, note that one passes from
$(A_1)$ to $(A_2)$ by replacing the rectangle defined by $\tau$ by 
the parallelogram defined by $\tau+1$, hence by {\bf 2.8}
and {\bf 2.9}, $\mu$ by $1-\mu$. From the geometric point of view this
is just one of the cases of Corollary {\bf 3.7} and can be described 
as applying half twists along $\gamma_1$ and $\gamma_2$.

To pass from $(A_1)$ to $(A_3)$ we 
replace $\tau$ by $\tau/(1-\tau)$, and hence $\mu$ by $1/\mu$. From
there the other transformations needed are clear, one passes from
$(A_2)$ to $(A_6)$ (resp. from $(A_3)$ to $(A_5)$) by a quarter
Dehn-twist along the median of the vertical cylinder (resp. horizontal
cylinder) and from $(A_4)$ to $(A_5)$ by two vertical half-twists.

The different values one obtains are
best expressed in terms of $\nu=a^2+1/a^2$ which is a modular
invariant (see [Si2] section {\bf 3} ). The values one obtains are
given in figure 15. 
\bigskip
\bigskip

\centerline{
 \SetLabels
\L(.11*1.05) $(A_1)$\\
\L(.43*1.05) $(A_2)$\\
\L(.795*1.05) $(A_3)$\\
\L(.11*.48) $(A_4)$\\
\L(.32*.2) $(A_5)$\\
\L(.795*.62) $(A_6)$\\
\L(.13*.63) $\nu$\\
\L(.4*.3) $\dfrac{2(6-\nu)}{2+\nu}$\\
\L(.805*.78) $-\nu$\\
\L(.08*-.12) $\dfrac{2(\nu+6)}{\nu-2}$\\
\L(.42*-.12) $\dfrac{2(\nu-6)}{\nu+2}$\\
\L(.77*-.12) $\dfrac{2(6+\nu)}{2-\nu}$\\
\endSetLabels
\AffixLabels{\psfig{figure=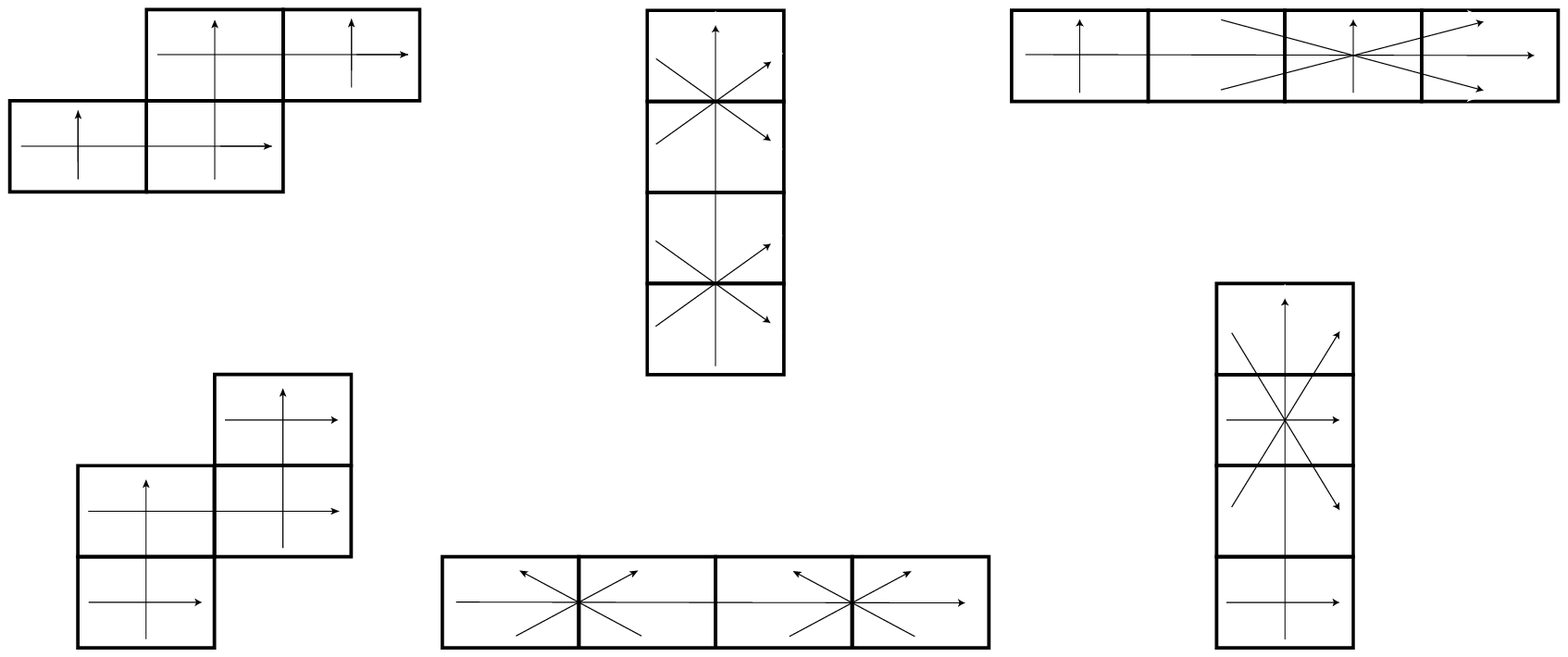,height=5cm}}}
\bigskip
\bigskip
\bigskip

\centerline{\bf Fig. 15}
\bigskip

\noindent{\bf Case $(B)$}. This family is Escb$_2(2)$ and we have an 
equation of the form $y^2=x\,(x^4+a\,x^2+1)$ (see section {\bf 4})
with
$$a=\frac{2\mu-4}{\mu}\ .\tag 5.13$$

There are three surfaces $(B_1)$, $(B_2)$ and $(B_3)$ in the
$\SL_2(\Z)$ orbit (see figure 16). One passes from $(B_1)$ to $(B_2)$ 
by two half-Dehn twists along the medians of the horizontal cylinders
or as above by replacing $\mu$ by $1-\mu$. In the same way one passes 
from $(B_1)$ to $(B_3)$ by replacing $\mu$ by $1/\mu$. This yields the 
different values indicated in figure 16.
\bigskip

\centerline{
 \SetLabels
\L(.11*.9) $(B_1)$\\
\L(.795*.8) $(B_3)$\\
\L(.315*.9) $(B_2)$\\
\L(.13*.08) $a$\\
\L(.39*-.28) $\dfrac{2(6-a)}{a+2}$\\
\L(.77*.0) $\dfrac{2(6+a)}{a-2}$\\
\endSetLabels
\AffixLabels{\psfig{figure=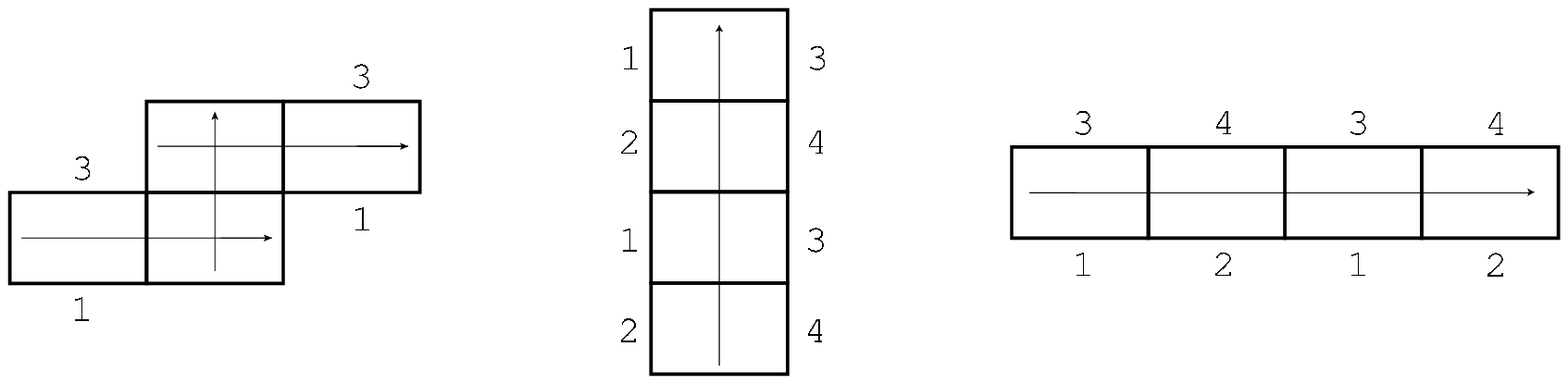,height=2.5cm}}}
\bigskip
\bigskip
\bigskip
\medskip

\centerline{\bf Fig. 16}
\bigskip

There are of course obvious similarities between case $(A)$ and case
$(B)$. These are explored in detail in [Si2] section {\bf 3}.
\bigskip

\noindent{\bf 5.14 Remark.} Here again we can be far more general. For
exactly the same reasons, the transformations 
$$a\mapsto\dfrac{2(6-a)}{a+2}\ \text{ and }\
a\mapsto\dfrac{2(6+a)}{a-2}$$
correspond to replacing in any of the escalator families  
$\tau$ by $1+\tau$ or $1/\tau$ respectively. These in turn correspond
to applying half-Dehn twists along the horizontal cylinders or the
vertical cylinders.
\bigskip

\noindent{\bf Case $(C)$}. 
We again label the rectangles in figure 13, $R_1$ to $R_4$, starting 
with the upper left.
We normalize so that the $x$-coordinate of the vertices of the rectangles
are the points at infinity. We normalize further so that the center of 
rectangle $R_1$ has $x$-coordinate 1 and the midpoint of the horizontal edge 
between rectangles $R_2$ and $R_3$ has $x$-coordinate 0. Note that since this
is not a Weierstrass point the midpoint of the lower edge of rectangle
$R_2$ (arc labeled 2 in figure 13) will also have $x$-coordinate 0. We will
call $a$, $b$ and $c$ the $x$-coordinate of, respectively, the midpoint
of the horizontal edges of rectangle $R_1$, the center of rectangle $R_2$ and
the midpoint of the upper horizontal edge of rectangle $R_3$ (arc labeled 3 in
figure 13). This is also the lower edge of rectangle $R_4$. 
Finally we will call $d_1$ and $d_2$ the $x$-coordinates of the
midpoints of the vertical edges of rectangles $R_3$ and $R_4$. 
The Weierstrass points are the points with $x$-coordinates 1, $a$,
$b$, $c$, $d_1$ and $d_2$.

With this notation we can choose for the map $f:\Bbb P^1\to \Bbb P^1$ that 
induces the covering map from the surface to the genus 1 curve obtained from  
one rectangle, the map
$$f(x)=\frac{x^2(x-a)(x-c)}{(1-a)(1-c)}\ .\tag 5.15$$
This map sends $1$, $b$, and two points with the same $x$-coordinate, to $1$.
Again this means that the equation $f(x)=1$ must have a double root 
(different from 1). This imposes conditions on $a$ and $c$, namely 
$$a=\frac{-(3\,t^2+4\,t+2)(t+2)}{t\,(t^2+2\,t+2)}\qquad
c=\frac{-(t^2+4\,t+6)(t+1)}{t^2+2\,t+2}\ .\tag 5.16$$
From this one recovers
$$b=-\frac{t^2+3\,t+2}t\ .\tag 5.17$$

Let $\lambda=f(d_i)$. In addition to the $d_i$ we again have two points 
with the same $x$-coordinate. Hence again $f(x)=\lambda$ must have a double 
root (with of course $\lambda\neq 1$). This imposes
$$\lambda=-\frac{(t^2-2)^2(3\,t^2+4\,t+2)^3(t^2+4\,t+6)^3}
{1024\,t^3(t^2+2\,t+2)^2(t+2)^3(t+1)^3}\ .\tag 5.18$$
This defines the genus 1 quotient and from  this on can easily compute 
the $d_i$. We find 
$$\multline
d_1=\frac{(t^2-2+(t^2+4\,t+2)\,i\sqrt{2})(t^2-2)}{4\,t\,(t^2+2\,t+2)},
\\  
d_2=\frac{(t^2-2-(t^2+4\,t+2)\,i\sqrt{2})(t^2-2)}{4\,t\,(t^2+2\,t+2)}
\ .\endmultline\tag 5.19$$

To obtain a full map we must introduce the points $p$ and $q$ which 
are mapped under $f$ of (5.15), to $f(d_1)=f(d_2)$ and 1 respectively.
We have 
$$p=-\frac{(t^2+4\,t+6)(3\,t^2+4\,t+2)}{4\,t\,(t^2+2\,t+2)}\ 
\text{ and }\ q=-\frac{2(t^2+3\,t+2)}{t^2+2\,t+2}\ .$$
With this a full map to the genus 1 curve is
$$(x,y)\mapsto\left(\frac{x^2(x-a)(x-c)}{(1-a)(1-c)},
y\frac{x\,(x-p)(x-q)}{\sqrt{(1-a)^3(1-c)^3}}\right)\ .$$  

This does not comply with our convention of section {\bf 2} but to
recover $\mu$ in terms of $\lambda$ we only need to set 
$$\mu=\frac{\lambda}{\lambda-1}\ .\tag 5.20$$
We will use this later. 

The solution we have found is of course far from optimal since solving (5.18)
for a specific value of $\lambda$ yields in general 16 solutions in $t$. 
On the other hand these solutions come in groups of 4. Namely if
$t$ is a solution then so are,
$$t,\ \frac2t,\ -\frac{t+2}{t+1},\ -\frac{2(t+1)}{t+2}\ .\tag 5.21$$
Moreover these 4 solutions yield isomorphic curves, since 
replacing $t$ by $\frac2t$ leaves $b$ fixed, exchanges $a$ and $c$ and 
exchanges the $d_i$, while replacing $t$ by $-\frac{t+2}{t+1}$ replaces 
$b$ by $1/b$, $a$ by $c/b$, $c$ by $a/b$ and the $d_i$ by $d_i/b$. In other 
words these depend on choices in the above computations.

We can use relation (5.21) to simplify the equation (5.18).
If we let $w$ be a root of
$$f_{\lambda}=(\lambda-1)x^4-(6\lambda+2)x^3+12\lambda\,x^2-(8\lambda-2)x+1\
,$$
or equivalently
$$f_{\mu}=x^4+(2-8\,\mu)x^3+12\mu\,x^2-(2+6\,\mu)x+\mu-1\ ,\tag 5.22$$
then solutions of the form (5.21) will be  
$$\frac{\sqrt{2}(u+1+\sqrt{2})}{u-1-\sqrt{2}},\ \text{ where }
\left(\frac{u^2+1}{u^2-1}\right)^2=w\ .\tag 5.23$$

On the other hand we cannot improve further. The reason is that for the
four surfaces in the $\SL_2(\Bbb Z)$-orbit we have the same repartition 
of Weierstrass points. Namely, with the same convention as before, it is
$(0,2,2,2)$.

To differentiate the cases we must take a closer look at equation (5.22).
If $\mu$ is real $>1$ then, $f_{\mu}=0$ will have two real roots,  
$w_1>2$ and $0<w_2<1/2$, and two complex conjugate roots $w_3$ and $w_4$. 
Computing $a$, $b$, $c$, $d_1$ and $d_2$ in terms of these roots we find 
that $a$, $b$ and $c$ will be real for $w_1$ while $d_1$ and $d_2$ are
complex conjugate, for $w_2$ on the other hand $b$, $d_1$ and $d_2$ will 
be real and $a$ and $c$ complex conjugate. Comparing with the real 
structures induced by reflection along the horizontal axis or the vertical 
axis in cases $(C_1)$ and $(C_4)$ we conclude that $w_1$ corresponds to 
$(C_1)$ and $w_2$ corresponds to $(C_4)$ (recall that $d_1$ and $d_2$ are 
on vertical edges while $a$ and $c$ are on horizontal edges).

In the absence of additional information we can not distinguish between 
$(C_2)$ and $(C_3)$ which are mirror images of each other and correspond 
to the two complex conjugate roots $w_3$ and $w_4$. 

Finally we note that we pass from $(C_1)$ to $(C_2)$ (resp.~ $(C_3)$)
by a third of a Dehn twist (resp.~ minus a third of a twist) along the 
geodesic median of the vertical cylinder of width 3 (and a full twist 
along the cylinder of width 1). We pass from $(C_2)$ to $(C_4)$ by a 
third of a Dehn-twist along the geodesic median of the horizontal 
cylinder of width 3.
\bigskip
\bigskip
\bigskip

\centerline{
 \SetLabels
\L(.1*1.2) $(C_1)$\\
\L(.34*1.2) $(C_2)$\\
\L(.605*1.2) $(C_3)$\\
\L(.845*1.2) $(C_4)$\\
\L(.085*-.2) $w_1>2$\\
\L(.355*-.2) $w_3$\\
\L(.59*-.2) $w_4=\overline{w}_3$\\
\L(.82*-.2) $0<w_2<\frac12$\\
\endSetLabels
 \AffixLabels{\psfig{figure=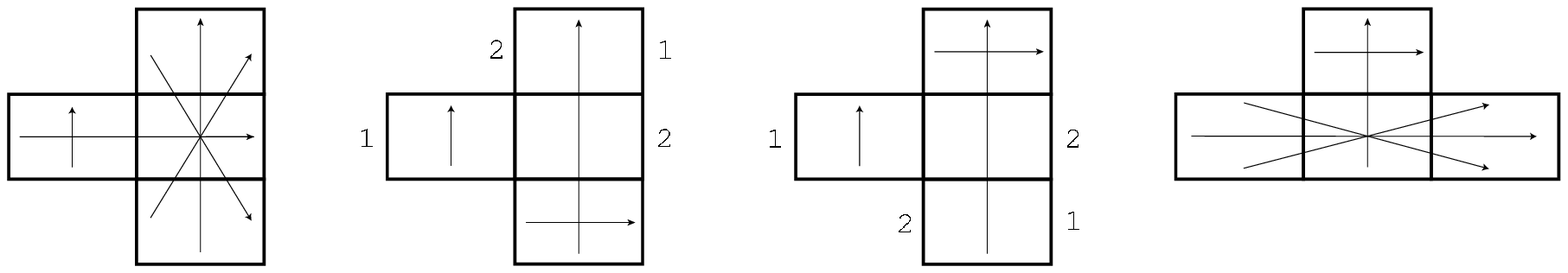,height=2.1cm}}}
\bigskip
\bigskip

\centerline{\bf Fig. 17}
\bigskip

\noindent{\bf Case $(D)$.} Label as before the rectangles in figure 13 
$R_1$ to $R_4$ starting with the lower left. The surface in this case 
has two non-hyperelliptic involutions with centers at the midpoints of 
the vertical edges of rectangles $R_1$ and $R_2$ for the first and at
the centers of the rectangles $R_3$ and $R_4$ for the second.

We normalize so that the $x$-coordinates of the vertices are $\pm i$, 
the centers of rectangles $R_3$ and $R_4$ are the points at infinity 
and the midpoints of the vertical edges of rectangles $R_1$ and $R_2$ 
have $x$-coordinate 0. We denote $\pm a$ the centers of rectangles
$R_1$ and $R_2$ and $\pm b$ the midpoints of the vertical edges of 
rectangles $R_3$ and $R_4$. 

These choices yield an equation of the curve in the form
$$y^2=(x^2-a^2)(x^2+1)(x^2-b^2)\ .\tag 5.24$$

A map from the curve to the genus 1 curve defined by one rectangle is
induced by 
$$f(x)=\frac{x^2(x^2-a^2)}{(x^2+1)^2}\ .\tag 5.25$$

Since by construction we have $f(b)=f(-b)=1$ we find
$$b=\frac{\pm i}{\sqrt{a^2+2}}\ .\tag 5.26$$

To find the complete equation of the genus 1 curve, we look for double 
roots of the equation $f(x)=\lambda$. This yields
$$\lambda=1-\mu=-\frac{a^4}{4\,(a^2+1)}\ ,\tag 5.27$$
where $\mu$ conforms with the convention of section {\bf 2}.

There are six different configurations in an $\SL_2(\Z)$ orbit (see
figure 18) and to differentiate them we again look at the possible
real structures and repartition of Weierstrass points among vertices, 
centers and horizontal and vertical edges. These are $(2,2,0,2)$ for
$D_1$ and $D_3$, $(2,0,2,2)$ for $D_2$ and $D_3$ and $(2,2,2,0)$ for
$D_5$ and $D_6$ (see figure 18).
\bigskip
\bigskip

\centerline{
 \SetLabels
\L(.10*1.05) $(D_1)$\\
\L(.36*1.05) $(D_2)$\\
\L(.67*1.05) $(D_3)$\\
\L(.9*1.05) $(D_4)$\\
\L(.24*.33) $(D_5)$\\
\L(.73*.33) $(D_6)$\\
\endSetLabels
\AffixLabels{\psfig{figure=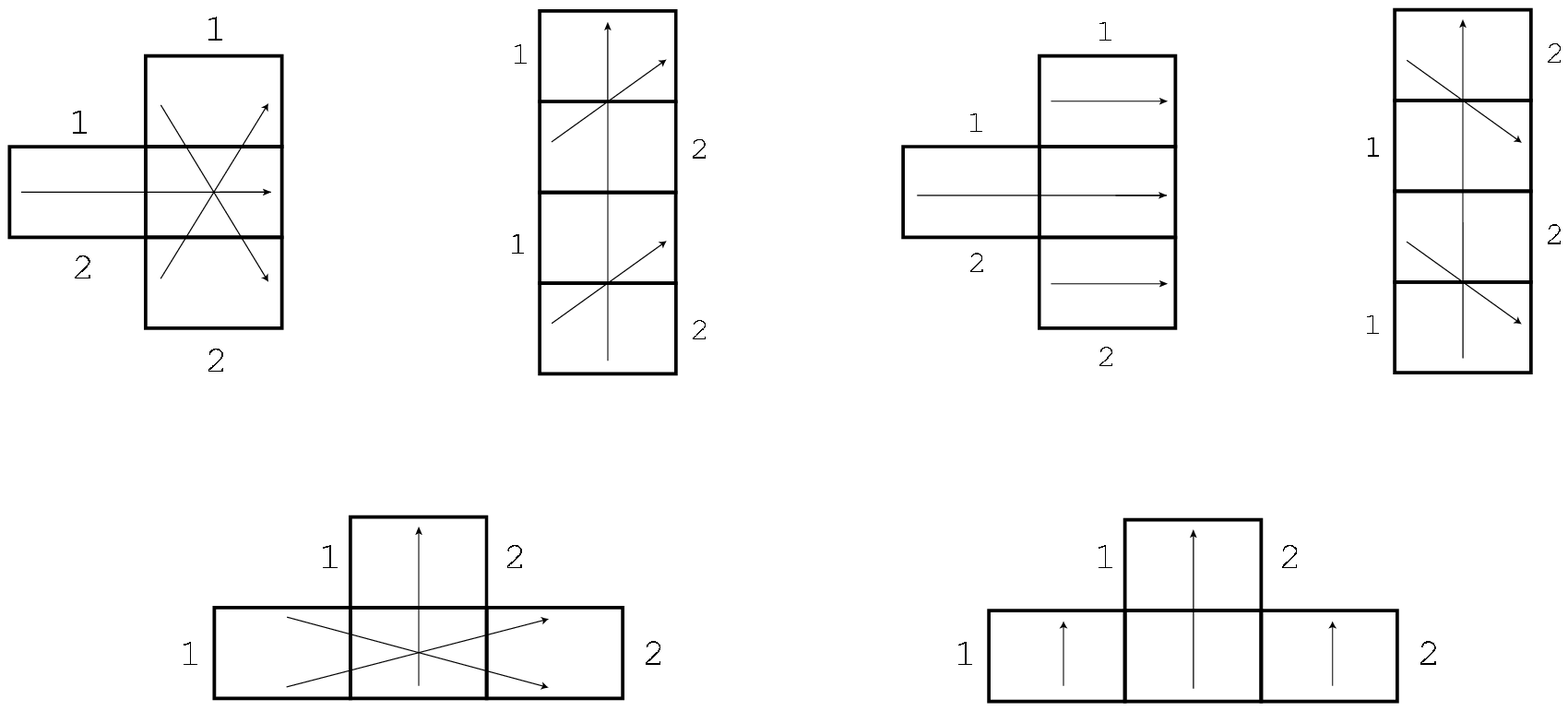,height=5cm}}}
\bigskip

\centerline{\bf Fig. 18}
\bigskip

Hence the surfaces in this family come in three pairs. In the first
$(D_1),(D_3)$ the two surfaces are real, but the first has only one
real component, while the second has three (for the real structure
induced by reflection along the vertical median of the cylinder of
width 4). For $\mu>1$ equation (5.27) has 2 real roots and two pure
imaginary roots
$$\pm\sqrt{2\,\mu-2+2\sqrt{\mu^2-\mu}}\ \ \text{ and }\ 
\pm\sqrt{2\,\mu-2-2\sqrt{\mu^2-\mu}}\ .\tag 5.28$$

Because of (5.24) and (5.26) this implies that the first corresponds to
$(D_1)$ while the second corresponds to $(D_3)$. The relation between 
$(D_5)$ and $(D_6)$ is similar and to find the corresponding value of
$a$ we only need, by Lemma {\bf 2.8}, to replace $\mu$ by $\mu/(\mu-1)$
in (5.28). The first value corresponds to $(D_5)$, the second to
$(D_6)$. 

To passes from $(D_1)$ to $(D_2)$ by a quarter Dehn-twist along the
geodesic median of the vertical cylinder of width 4 and similarly from 
$(D_2)$ to $(D_3)$, from $(D_3)$ to $(D_4)$ and from $(D_4)$ to
$(D_1)$. Hence the values of $a$ for $(D_2)$ and $(D_4)$ are obtained
by replacing $\mu$ by $1/\mu$ in (5.28) (see {\bf 2.8}).

Finally note that we pass from $(D_2)$ to $(D_5)$ by a quarter
Dehn-twist along the geodesic median of the horizontal cylinder of
width 4.
\medskip

We note here that combining (5.24), (5.25) and (5.28), we obtain
equations

$$\aligned
y^2&=(x^2+1)\left(x^2-2\,\left(\mu-1+\sqrt{\mu^2-\mu}\right)\right)
\left(x^2+\frac{\mu-\sqrt{\mu^2-\mu}}{2\,\mu}\right)\ 
\text{ for }\ (D_1)\\
y^2&=(x^2+1)\left(x^2-2\,\left(\mu-1-\sqrt{\mu^2-\mu}\right)\right)
\left(x^2+\frac{\mu+\sqrt{\mu^2-\mu}}{2\,\mu}\right)\ 
\text{ for }\ (D_3)\endaligned\tag 5.29$$
with similar relations between $(D_2)$ and $(D_4)$ and between $(D_5)$ 
and $(D_6)$.

\bigskip

\subhead 6. Remarks on curves defined over number fields\endsubhead
\medskip

An immediate consequence of (5.29) is that if $\mu$ is a general
squarefree integer, then $(D_1)$ and $(D_3)$ are Galois conjugate
in $\Bbb Q[\sqrt{\mu}]$. The interesting point here is that one
passes from $(D_1)$ to $(D_3)$ by a fractional Dehn-twist. There seems
to be many more instances where such a phenomenon occurs. Simple
examples arise for a real quadratic number field $K=\Bbb Q[\sqrt{d}]$ 
if $\mu$ is a unit of norm 1 or if $\mu=\frac12+b\,\sqrt{d}$
since in the first case $\mu$ and $1/\mu$ are conjugate in
$K|\Bbb Q$ while in the second $\mu$ and $1-\mu$ are conjugate.
Examples of the first are $\mu=9+4\sqrt{5}$ and $\mu=97+56\sqrt{3}$
while $\mu=(2+\sqrt{5})/4$ and $\mu=(12+7\sqrt{3})/24$ are examples
for the second (see tables in Appendix).

The phenomenon is not restricted to quadratic extensions. Forgetting
the reference to $\mu$ the relation between three surfaces in 
the same $\PSL_2(\Z)$ orbit of case $(B)$ can be expressed as follows:
Let $a_1$, $a_2$ and $a_3$ be the three corresponding parameters of the
form (5.13), then the $a_i$ are the roots of an equation of the form   
$$x^3-\alpha\,x^2-36\,x+4\,\alpha\tag 6.1$$
and again one passes from one to the others by fractional Dehn-twists.

Similarly the six elements in a same $\PSL_2(\Z)$ orbit in case $(A)$
are obtained by letting the parameter $\nu$ be one of the roots of
 $$x^6-(\alpha+72)x^4+(8\,\alpha+1296)x^2-16\,\alpha\ .\tag 6.2$$

For case $(C)$, rewriting (5.24), the four surfaces are obtained by 
letting $w$ be one of the roots of
$$x^4+(2-8\,\alpha)x^3+12\alpha\,x^2-(2+6\,\alpha)x+\alpha-1\ .\tag
6.3$$

The situation is a little less satisfactory here since we were not
able to express directly equations for these surfaces in terms of the
parameter $w$. On the other hand one can show, although the
computations are too long to be presented here, that the coordinates
of the isomorphy classes of the surfaces in the Igusa moduli space 
(see [Ig]) can be expressed in terms of the parameter $w$.
\bigskip

\subhead Appendix: Values of $\mu$ for some equiquadrangles\endsubhead
\medskip

The values given in the tables that follow where computed numerically using 
a variant of the method described 
in [Bu-Si2] and should be considered as conjectured values. This being said
most of these values are in fact exact and can be
deduced from the known exact uniformization of certain curves (see for
example [Ai-Si], [Bu-Si1], [Bu-Si2], [Si1] or [Si2]). In fact these exact 
cases often come from surfaces at the intersection of two families. For example
the curve defined by $y^2=x^6+1$ is in Escb$_1(2)$ and obtained from squares. 
But it also has an automorphism of order 4 and hence is also in Escb$_2(2)$, 
the corresponding equation is $y^2=x\,(x^4-\frac{10}3\,x^2+1)$ and from this
we obtain for angle $\pi/4$ that $L=\sqrt{2}$ yields $\mu=\frac43$. In a 
similar fashion the curve with equation $y^2=x\,(x^4+1)$ is isomorphic to 
the one with equation $y^2=x^6-5\,i\,\sqrt{2}\,x^3+1$ and from this we can
deduce that the surface defined by the trace triple (see below) 
$(3+2\sqrt{2},\, 4+2\sqrt{2},\,4+2\sqrt{2})$ yields 
$\mu=\frac12-\frac{5\sqrt{2}}4\,i$. The computations are generally more 
involved but the method is the same. We give one last example: surfaces 
with Fenchel-Nielsen coordinates  $(\ell,\frac12,\ell,\frac12,\ell,\frac12)$ 
admit a second pants decomposition with coordinates 
$(\ell',\frac12,\ell,0,\ell',\frac12)$ (see [Si2]) where 
$L'=\cosh(\ell'/2)=\frac{3L-1}{2(L-1)}$ if $L=\cosh(\ell/2)$. We have
 $\ell'=\ell$ for $L=(5+\sqrt{17})/4$ and from this information it is 
possible to compute the $\mu$ for $(1+\sqrt{17})/4=\cosh(\ell/4)$
in the $\pi/3$ table.
\bigskip

In the first two tables $L$ denotes the hyperbolic cosine of the
hyperbolic half-length of the horizontal median of the equiquadrangle
with interior angle $\pi/n$ and $\mu$ is the invariant defined in
(2.7) or {\bf 2.9} 3. All the corresponding surfaces are with twist 
parameter zero.

We complete this list with a few examples with nonzero twist parameter.
These are best described in terms of a trace triple and we use the 
convention of [Ac-Na-Ro], that is we describe them with a triple 
$(x^2,y^2,z^2)$, where $x=\trace(A)$, $y=\trace(B)$ and 
$z=\trace(A\cdot B^{-1})$, $A$ and $B$ being the generators of the group.
For the surfaces with zero twist the corresponding trace triple is 
$(4L^2,\,4L'{}^2,\, 4L^2L'{}^2)$, $L'$ as in Lemma {\bf 2.1}.

It should be noted that the examples presented in this appendix all 
have arithmetic Fuchsian groups (see tables in [Ta] and [Ac-Na-Ro]).
 On the other hand for general $n$, and in particular large enough $n$, 
the groups described in {\bf (2.3)}, {\bf (2.4)} and {\bf (2.5)} are not 
arithmetic Fuchsian groups. They are however subgroups of triangle groups,
and this is in accordance with a conjecture of Chudnovsky and Chudnovsky 
([Ch-Ch] section 7) that if a curve defined over 
$\overline{\Bbb Q}$ has a Fuchsian group $G$ in 
$\PSL_2(\Bbb R\cap\overline{\Bbb Q})$ then $G$ is either an arithmetic 
group or a subgroup of a triangle group.
\bigskip
\newpage

\vbox{\tabskip=0pt\offinterlineskip
\def\tablerule{\noalign{\hrule}}
\def\yop{\vphantom{$\dfrac{\frac{\sqrt{1^2}^2}{2}}{\frac{1}{\sqrt{2}}}$}}
\halign to355pt{\strut#&\vrule#\tabskip=1em plus2em&\hfil#\hfil&\vrule#&
\hfil#\hfil&\vrule#&\ \vrule#&\hfil#\hfil&\vrule#&\hfil#\hfil&\vrule#
\tabskip=0pt\cr\tablerule
&&\multispan3 \hfil {\bf For angle $\pmb{\pi}/{\bold 3}$}\hfil&&&
\multispan3 \hfil {\bf For angle $\pmb{\pi}/{\bold 4}$}\hfil &\cr
\tablerule
&&\omit\hidewidth $L$\hidewidth &&\omit\hidewidth $\mu$\hidewidth &&
&\omit\hidewidth $L$\hidewidth &&\omit\hidewidth $\mu$\hidewidth&\cr
\tablerule
&&\yop $\sqrt{2+\sqrt{2}}$ && $\dfrac{837+1107\sqrt{2}}{2401}$ &&&
$\dfrac{\sqrt{10+2\sqrt{17}}}2$&& $\dfrac{1151-217\sqrt{17}}{256}$&\cr
\tablerule
&&\yop $\dfrac{\sqrt{7+\sqrt{17}}}2$&& $\dfrac{23+\sqrt{17}}{27}$&&&
$\dfrac{\sqrt{2}+\sqrt{6}}2$&& $\dfrac{12+7\sqrt{3}}{24}$&\cr
\tablerule
&&\yop $\dfrac{\sqrt{6+2\sqrt{3}}}2$ && $27-15\sqrt{3}$  &&& 
$\sqrt{3}$&& $\dfrac{128}{125}$&\cr
\tablerule
&&\yop$\sqrt{2}$ && $\dfrac{27}{25}$&&&  $\dfrac{1+\sqrt{5}}2$&& 
$\dfrac{2+\sqrt{5}}4$&\cr
\tablerule
&&\yop $\dfrac{\sqrt{5+\sqrt{5}}}2$&& $\dfrac{32}{27}$&&&
$\dfrac{\sqrt{5+\sqrt{17}}}2$ && $\dfrac{897-217\,\sqrt{17}}{2}$&\cr
\tablerule
&&\yop $\dfrac{\sqrt{4+2\sqrt{2}}}2$&& $\dfrac{1564+1107\sqrt{2}}{2401}$&&&
$\sqrt{2}$&& $\dfrac43$&\cr
\tablerule
&&\yop $\dfrac{1+\sqrt{17}}4$&& $\dfrac{621+27\sqrt{17}}{512}$&&&
$\dfrac{\sqrt{18+2\sqrt{33}}}4$ && $\dfrac{283+21\sqrt{33}}{256}$&\cr
\tablerule
&&\yop  $\dfrac{\sqrt{6}}2$&& $2$&&&$\dfrac{\sqrt{4+2\sqrt{2}}}2$&& $2$&\cr
\tablerule
&&\yop $\dfrac{\sqrt{14+2\sqrt{17}}}4$&& $-108+27\sqrt{17}$&&&
$\dfrac{\sqrt{3}+\sqrt{11}}4$&& $\dfrac{9+7\sqrt{33}}{18}$&\cr
\tablerule
&&\yop$\dfrac{\sqrt{4+\sqrt{2}}}2$&& $\dfrac{58+41\sqrt{2}}{27}$&&&
$\dfrac{\sqrt{6}}2$&& $4$&\cr
\tablerule
&&\yop $\dfrac{\sqrt{2}+\sqrt{10}}4$&& $\dfrac{32}5$&&&
$\dfrac{\sqrt{14+2\,\sqrt{17}}}4$&& $\dfrac{1151+217\,\sqrt{17}}{256}$&\cr
\tablerule
&&\yop$\dfrac{\sqrt{5}}2$&& $\dfrac{27}2$&&&
$\dfrac{\sqrt{2}+\sqrt{10}}4$&& $9+4\sqrt{5}$&\cr
\tablerule
&&\yop $\dfrac{\sqrt{3+\sqrt{3}}}2$ && $27+15\sqrt{3}$&&&
$\dfrac{\sqrt{5}}2$&& $\dfrac{128}3$&\cr
\tablerule
&&\yop$\dfrac{\sqrt{10+2\sqrt{17}}}4$&& $109+27\sqrt{17}$&&&
$\dfrac{\sqrt{3+\sqrt{3}}}2$&& $97+56\sqrt{3}$&\cr
\tablerule
&&\yop $\dfrac{\sqrt{3+\sqrt{2}}}2$ && $\dfrac{1566+1107\sqrt{2}}2$  
&&& $\dfrac{\sqrt{10+2\sqrt{17}}}4$&&
$\dfrac{897+217\sqrt{17}}2$&\cr
\tablerule
&&\multispan7  \hfil {\qquad\qquad\qquad\bf For angle $\pmb{\pi}/{\bold
5}$}\hfil &&\cr
\tablerule
&&\yop $ \dfrac{\sqrt{8+2\sqrt{5}}}2$&& $\dfrac{65+29\sqrt{5}}{125}$
&&& $\dfrac{\sqrt{2}+\sqrt{10}}4$ && $\dfrac{1621+725\sqrt{5}}{121}$&\cr
\tablerule
&&\yop $\dfrac{\sqrt{5+\sqrt{5}}}2$ && $2$ &&& \multispan3 &\cr
\tablerule\cr
}}
\newpage

\vbox{\tabskip=0pt\offinterlineskip
\def\tablerule{\noalign{\hrule}}
\def\yop{\vphantom{$\dfrac{\frac{\sqrt{1^2}^2}{2}}{\frac{1}{\sqrt{2}}}$}}
\halign to350pt{\strut#&\vrule#\tabskip=1em plus2em&\hfil#\hfil&\vrule#&
\hfil#\hfil&\vrule#&\ \vrule#&\hfil#\hfil&\vrule#&\hfil#\hfil&\vrule#
\tabskip=0pt\cr
\tablerule
&&\multispan3 \hfil {\bf For angle $\pmb{\pi}/{\bold 6}$}\hfil 
&&&\multispan3 \hfil {\bf For angle $\pmb{\pi}/{\bold 8}$}\hfil &\cr
\tablerule
&&\omit\hidewidth $L$\hidewidth &&\omit\hidewidth $\mu$\hidewidth &&
&\omit\hidewidth $L$\hidewidth &&\omit\hidewidth $\mu$\hidewidth&\cr
\tablerule
&&\yop$\dfrac{\sqrt{3}+\sqrt{7}}2$&& $128-48\sqrt{7}$&&&
$\sqrt{2+\sqrt{2}}$&& $\dfrac{-4+8\sqrt{2}}7$&\cr
\tablerule
&&\yop$2$&& $\dfrac{81}{80}$&&&
$\dfrac{\sqrt{6+2\sqrt{2}}}2$&& $\dfrac{3+2\sqrt{2}}4$ &\cr
\tablerule 
&&\yop $\dfrac{3\sqrt{2}+\sqrt{10}}4$ && $\dfrac{4096-1216\sqrt{10}}{243}$ &&&
$\sqrt{\cos(\pi/8)+1}$&& $2$&\cr
\tablerule 
&&\yop $ \dfrac{\sqrt{6+2\sqrt{7}}}2$ && $\dfrac{512-160\sqrt{7}}{81}$ &&&
$\dfrac{\sqrt{4+2\sqrt{2}}}2$&& $\dfrac{11+8\sqrt{2}}7$&\cr
\tablerule
&&\yop$\dfrac{\sqrt{10}}2$&& $\dfrac{32}{27}$&&&
$\dfrac{\sqrt{4+\sqrt{2}}}2$&& $12+8\sqrt{2}$&\cr
\tablerule
&&\yop $\sqrt{2}$&& $\dfrac{81}{49}$&&&\multispan3 \hfil {\bf For
angle ${\bold 0}$}\hfil &\cr
\tablerule
&&\yop $\dfrac{1+\sqrt{3}}2$&& $2$&&&
$\sqrt{5}$&& $\dfrac{125-55\sqrt{5}}2$&\cr
\tablerule
&&\yop $\dfrac{\sqrt{7}}2$&& $\dfrac{81}{32}$&&&
$\sqrt{3}$&& $\dfrac98$&\cr
\tablerule
&&\yop $\dfrac{\sqrt{6}}2$&& $\dfrac{32}{5}$&&&
$\sqrt{2}$&& $2$&\cr
\tablerule
&&\yop$\dfrac{\sqrt{3+\sqrt{7}}}2$&& $\dfrac{512+160\sqrt{7}}{81}$
&&& $\dfrac{\sqrt{6}}2$&& $9$&\cr
\tablerule
&&\yop$\dfrac{\sqrt{2}+\sqrt{10}}4$&& $\dfrac{4096+1216\sqrt{10}}{243}$&&&
$\dfrac{\sqrt{5}}2$&& $\dfrac{125+55\sqrt{5}}2$ &\cr
\tablerule
&&\yop$\dfrac{\sqrt{5}}2$&& $81$&&& \multispan3&\cr
\tablerule
&&\yop$\dfrac{\sqrt{3}+\sqrt{7}}4$&& $128+48\sqrt{7}$&&&
\multispan3 &\cr
\tablerule\cr
}}
\bigskip

\noindent For surfaces with a nonzero twist parameter, in terms of the trace
triple $(x^2,y^2,z^2)$.
\bigskip

\noindent For $\pi/3$\ \ (i.e.~ sum of the interior angles equal to $4\pi/3$).

$(3+\sqrt{7},\,4+\sqrt{7},\,5+\sqrt{7})$, \ 
$\mu=\dfrac{59+17\sqrt{7}}{54}(1+i)$.

$(3+2\sqrt{2},\, 4+2\sqrt{2},\,4+2\sqrt{2})$, \ 
$\mu=\dfrac12-\dfrac{5\sqrt{2}}4\,i$.
\smallskip

$(9,6,6)$, \ $\mu=1/2$
\smallskip

$(1+4\rho+4\rho^2,\,1+4\rho+4\rho^2,\,1+4\rho+4\rho^2)$, $\rho=\cos(2\pi/9)$, 
\ $\mu=(1-\sqrt{3}\,i)/2$.

\bigskip

\noindent For $\pi/4$ \ \ (i.e.~ sum of the interior angles equal to $\pi$).

$(9,7,7)$,\ $\mu=\dfrac12-\dfrac{13\sqrt{7}}{98}\, i$.

$(4+2\sqrt{3},\, 4+2\sqrt{3},\, 4+2\sqrt{3})$, \  $\mu=(1-\sqrt{3}\,i)/2$.
\smallskip

$(3/2+\sqrt{2},\, 4+2\sqrt{2},\, 4+2\sqrt{2})$, \ $\mu=1/2$.
\smallskip

$(3+\sqrt{6},\,5+2\sqrt{6},\,6+2\sqrt{6})$, \ 
$\mu=(5+2\sqrt{6})(1-\sqrt{3}\,i)/2$.
\bigskip

\noindent For $\pi/6$\ \ (i.e.~ sum of the interior angles equal to $2\pi/3$).

$(9,8,8)$,\ $\mu=\dfrac12-\dfrac{7\sqrt{2}}{16}\,i$.

$(7/2+3\sqrt{5}/2,\, 5+2\sqrt{5},\,5+2\sqrt{5})$, \
$\mu=\dfrac12-\dfrac{79\sqrt{5}}{80}\,i$.

$(3+2\sqrt{2},\,6+4\sqrt{2},\,6+4\sqrt{2})$, \ 
$\mu=\dfrac12-\dfrac{11\sqrt{2}}{11}\,i$.

$(7+4\sqrt{3},\,4+2\sqrt{3},\,4+2\sqrt{3})$, \ $\mu=1/2$.
\smallskip 

$(\rho^2,\rho^2,\rho^2)$, $\rho=1+2\cos(\pi/9)$, \ $\mu=(1-\sqrt{3}\,i)/2$.

\end

\end

\end

 Because of this it does not make
much sense to use Fenchel-Nielsen coordinates and it is better to directly 
give the generators of the Fuchsian group. With the notations of section 
{\bf ?} these will be 
$$\langle A^2,B,AB^3A,(AB)A(AB)^{-1},(AB^2)A(AB^2)^{-1}\rangle\ .
\tag 4.22$$

\end
 
we note that we can always choose for $E$ an
equation of the form $y^2=x^4+a\,x^2+1$. This curve is isomorphic to the one 
dfined by
$$v^2=u\,(u-1)\,(u-\mu)\ \text{ where }\ \mu=\frac4{2-a}\tag 4.3$$
an explicit isomorphism being
$$x=\frac{v\,\sqrt{\mu^3}}{\mu\,u\,(u-\mu)},\ 
y=\frac{u^2-2\,u+\mu}{u\,(u-\mu)}\ .$$

\bye

%% file: labelfig.tex
  \ifx \LabelFigloaded\MYundefined
  \else
    \immediate\write16{ !!! LabelFig.tex ALREADY loaded.}
   \fi

  \def\LabelFigloaded{\relax}


  \chardef\LabelFigCatAt\the\catcode`\@
  \catcode`\@=11

 \let\LabelFigwlog@ld\wlog
 \def\wlog#1{\relax}

 \ifx\\\MYundefined@
    \let\\\relax
 \fi

 \def\N@wif{\csname newif\endcsname }
 \def\Temp@ {\N@wif\ifIN@}
 \ifx\INN@\MYundefined@
    \else \let\Temp@\relax
 \fi
 \Temp@

  \def\IN@{\expandafter\INN@\expandafter}
  \long\def\INN@0#1@#2@{\long\def\NI@##1#1##2##3\ENDNI@
    {\ifx\m@rker##2\IN@false\else\IN@true\fi}%
     \expandafter\NI@#2@@#1\m@rker\ENDNI@}
  \def\m@rker{\m@@rker}
 \def\Shifted@@#1#2#3{\setbox0=\hbox{#3}%
   \raise -\dp0\vbox {\kern-#2%
       \hbox {\kern#1\box0\kern-#1}%
           \kern#2}}

 \newcount\gridcount
 \newbox\auxGridbox@ \newbox\hGridbox@ \newbox\vGridbox@
 \newbox\Labelbox@ \newbox\auxLabelbox@
 \newbox\Coordinatebox@
 \newtoks\Labeltoks@
 \newdimen\Wdd@ \newdimen\Htt@

 \def\hRule@{\advance\gridcount -2%
   \vskip-.2pt\hrule\vskip-.2pt\vfil
   \llap{\smash{\raise -2.5pt
     \hbox{.\number\gridcount\kern2pt}}}%
           \vskip-.2pt\hrule\vskip-.2pt\vfil}

 \def\vRule@{\advance\gridcount 2%
   \hskip-.2pt\vrule\hskip-.2pt\hfil
   \setbox\auxGridbox@=\vbox to 0pt
      {\vskip \Htt@\vskip 2pt
        \hbox{\kern-3.5pt.\number\gridcount}\vss}%
      \wd\auxGridbox@=0pt \box\auxGridbox@
   \hskip-.2pt\vrule\hskip-.2pt\hfil}

 \def\PlaceGrid@@{\gridcount=10%
  \setbox\hGridbox@=%
    \hbox{\hbox{\GridSpider@{\hskip-.4pt\vrule}%
             \vbox to \Htt@{\offinterlineskip\parindent=\z@%
                \GridSpider@{\vskip-.4pt\hrule}\vfil
                \hRule@\hRule@\hRule@\hRule@
                  \vskip-.2pt\hrule\vskip-.2pt\vfil
                \hbox to \Wdd@{\hfil}%
             \GridSpider@{\hrule\vskip-.4pt}}%
         \GridSpider@{\vrule\hskip-.4pt}}}%
  \gridcount=0%
  \setbox\vGridbox@=
     \hbox{\vbox{\offinterlineskip\parindent=0pt\hsize=0pt
       \GridSpider@{\vskip-.4pt\hrule}%
             \hbox to \Wdd@{%
                \GridSpider@{\hskip-.4pt\vrule}\hfil
                \vtop to \Htt@{\vfil}%
                     \vRule@\vRule@\vRule@\vRule@
                     \hskip-.2pt\vrule\hskip-.2pt\hfil
             \GridSpider@{\vrule\hskip-.4pt}}%
         \GridSpider@{\hrule\vskip-.4pt}}}%
  \wd\hGridbox@=0pt\ht\hGridbox@=0pt
  \wd\vGridbox@=0pt\ht\vGridbox@=0pt
 \hbox{\box\hGridbox@\box\vGridbox@}%
  }

 \def\SetLabels#1\endSetLabels{%
   \Labeltoks@={#1}}

 \def\GridSpider@#1{#1}
 \let\PlaceGrid@\relax
 \def\ShowGrid{\let\PlaceGrid@\PlaceGrid@@}

 \def\bAdjust@@{%
     \setbox\auxLabelbox@=\hbox{\raise \dp\auxLabelbox@
            \box\auxLabelbox@}}
 \def\bAdjust@{\let\vAdjust@\bAdjust@@}

 \def\tAdjust@@{%
     \setbox\auxLabelbox@=\hbox{\raise-\ht\auxLabelbox@
            \box\auxLabelbox@}}
 \def\tAdjust@{\let\vAdjust@\tAdjust@@}

 \let\vAdjust@\relax

 \def\lAdjust@{\let\hAdjust@\rlap}
 \def\rAdjust@{\let\hAdjust@\llap}

 \let\hAdjust@\relax\let\vAdjust@\relax

 \def\FetchLabel@#1(#2*#3)#4\\#5\endFetchLabel@{%
     \ignorespaces#1\unskip
     \Labeltoks@={#5}%
     \setbox\auxLabelbox@=%
        \hbox to 0pt{\hss\ignorespaces\hAdjust@
          {\ignorespaces#4\unskip}\hss}%
     \vAdjust@
     \let\hAdjust@\relax\let\vAdjust@\relax
     \setbox\Labelbox@=\hbox to 0pt{%
       \box\Labelbox@
       \Shifted@@{#2\Wdd@}{#3\Htt@}{\box\auxLabelbox@}}%
     \ht\Labelbox@=0pt\dp\Labelbox@=0pt
     }

 \def\PlaceLabels@@{\bgroup\def\Cr@{\\}%
     \let\L\lAdjust@\let\R\rAdjust@
     \let\B\bAdjust@\let\T\tAdjust@
     \loop
     \IN@0\Cr@ @\the\Labeltoks@ @\relax
     \ifIN@ \expandafter
       \FetchLabel@\the\Labeltoks@\endFetchLabel@
     \repeat
     \box\Labelbox@\egroup}%

 \let \PlaceLabels@\PlaceLabels@@

 \def\AffixLabels#1{\setbox\Coordinatebox@=\hbox{#1}%
      \Wdd@=\wd\Coordinatebox@ \Htt@=\ht\Coordinatebox@
      \advance\Htt@ \dp \Coordinatebox@
      \hbox{\copy\Coordinatebox@\kern-\Wdd@%
           \Shifted@@{0pt}{-\dp\Coordinatebox@}%
            {\PlaceGrid@\PlaceLabels@}%
           \kern\Wdd@}}

   \let\wlog\LabelFigwlog@ld 
   \catcode`\@=\LabelFigCatAt  


%% file: multigeodesic.bbl
\bigskip


\Refs
\refstyle{A}
\widestnumber\key{Ac-Na-Ro}
\ref\key Ac-Na-Ro
\by P. Ackermann, M. N\"a\"at\"anen and G. Rosenberger
\paper The arithmetic Fuchsian groups with signature $(0;2,2,2,q)$
\jour Research and Exposition in Math. \vol 27 \yr 2003 \pages 1--9
\endref
\ref\key Ai-Si
\by A. Aigon and R. Silhol
\paper Hyperbolic hexagons and algebraic curves in genus 3 
\jour J. London Math. Soc. \yr 2002 \vol 66 \pages 671--690 
\endref
\ref\key{Bea}
 \by A. Beardon
\book The Geometry of Discrete Groups
\publ Springer G.T.M. 91\publaddr Berlin Heidelberg New York
\yr 1991
\endref
\ref\key Bu
\by P. Buser
\book Geometry and Spectra of Compact Riemann Surfaces
\publ Birkh\"auser\publaddr Boston Basel Berlin
\yr 1992
\endref
\ref\key{Bu-Se}
\by P. Buser and M. Sepp\"al\"a
\paper Real structures of Teichm\"uller spaces, Dehn twists, and
moduli spaces of real curves
\jour Math. Z. \vol 232 \yr 1999 \pages 547--558
\endref
\ref\key{Bu-Si1}
\by P. Buser and R. Silhol
\paper Geodesics, periods and equations of real hyperelliptic curves
\jour Duke Math. J. 
\vol 108 \yr 2001 \pages 211--250
\endref
\ref\key{Bu-Si2}
\by P. Buser and R. Silhol
\paper Some remarks on the uniformizing function in genus 2
\jour Geometria Dedicata \yr 2005 \vol 115 \pages 121--133
\endref
\ref\key Ch-Ch
\by D.V. Chudnovsky and G.V. Chudnovsky
\paper Computer algebra in the service of mathematical physics and number 
theory \jour Computers and Mathematics, Proceedings Stanford 1986, Dekker 
1990 \pages 109--232
\endref
\ref\key{Hu-Le}
\by P. Hubert and S. Leli\`evre
\paper Prime arithmetic Teichm\"uller discs in $\Cal H(2)$
\jour Israel J. of Math. \vol 151 \yr 2006 \pages 281--321
\endref
\ref\key{Ig} 
\by J. I. Igusa
\paper Arithmetic variety of moduli for genus two \jour Ann. of Math. 
\vol 72 \yr 1960 \pages 612--649
\endref
\ref\key{Ku-N\"a}
\by T. Kuusalo and M. N\"a\"at\"anen
\paper Geometric uniformization in genus 2
\jour Annales Acad. Sci. Fennic\ae \vol 20 \yr 1995 \pages 401--418
\endref
\ref\key Mc
\by C. T. McMullen
\paper Billiards and Teichm\"uller curves on Hilbert modular surfaces
\jour J. Amer. Math. Soc. \yr 2003 \vol 16 \pages 857--885
\endref
\ref\key M\"o
\by M. M\"oller \paper Teichm\"uller curves, Galois actions and
$\widehat{GT}$-relations \jour Math. Nachr. \vol 278 \yr 2005 \pages
1061--1077 
\endref
\ref\key Ne
\by Z. Nehari
\book Conformal Mapping
\publ McGraw-Hill\publaddr New York Toronto London
\yr 1952
\endref
\ref\key Sch
\by G. Schmith\"usen
\paper Examples for Veech groups of origamis
\jour  Contemporary Math. \vol397 \yr 2006 \pages 193--206
\endref
\ref\key Si1
\by R. Silhol
\paper Hyperbolic Lego and equations of algebraic curves 
\jour Contemporary Math. \vol 311 \yr 2002 \pages 313--334
\endref
\ref\key Si2
\by R. Silhol
\paper
On some one parameter families of genus 2 algebraic curves and half twists
\jour to appear in Commentarii Math. Helvetici\yr 2007 \vol 82
\endref
\ref\key Ta
\by K. Takeuchi \paper Arithmetic Fuchsian groups with signature $(1;e)$
\jour J. Math. Soc. Japan \vol 35 \yr1983 \pages 381--407 
\endref
